\documentclass{article}

\usepackage{PRIMEarxiv}

\usepackage[utf8]{inputenc} % allow utf-8 input
\usepackage[T1]{fontenc}    % use 8-bit T1 fonts
\usepackage{hyperref}       % hyperlinks
\usepackage{url}            % simple URL typesetting
\usepackage{booktabs}       % professional-quality tables
\usepackage{amsfonts}       % blackboard math symbols
\usepackage{nicefrac}       % compact symbols for 1/2, etc.
\usepackage{microtype}      % microtypography
\usepackage{lipsum}
\usepackage{fancyhdr}       % header
\usepackage{graphicx}       % graphics
\graphicspath{{media/}}     % organize your images and other figures under media/ folder

\usepackage{amssymb}
\usepackage{amsmath,amsthm}

\newtheorem{theorem}{Theorem}[section]

\newtheorem{proposition}{Proposition}[section]
\newtheorem{lemma}{Lemma}[section]
\newtheorem{remark}{Remark}
\newtheorem{corollary}{Corollary}[section]
\newtheorem{example}{Example}[section]

\title{The Optimal Strategy against Hedge Algorithm in Repeated Games
}

\author{
  Xinxiang Guo, Yifen Mu \\
  The Key Laboratory of Systems and Control, Academy of Mathematics and Systems Science \\
  Chinese Academy of Sciences \\
  Beijing, China\\
  \texttt{\{guoxinxiang, mu\}@amss.ac.cn} \\
}

\begin{document}
\maketitle

\begin{abstract}
This paper aims to solve the optimal strategy against a well-known adaptive algorithm, the Hedge algorithm, in a finitely repeated $2\times 2$ zero-sum game. In the literature, related theoretical results are very rare. To this end, we make the evolution analysis for the resulting dynamical game system and build the action recurrence relation based on the Bellman optimality equation. First, we define the state and the State Transition Triangle Graph (STTG); then, we prove that the game system will behave in a periodic-like way when the opponent adopts the myopic best response. Further, based on the myopic path and the recurrence relation between the optimal actions at time-adjacent states, we can solve the optimal strategy of the opponent, which is proved to be periodic on the time interval truncated by a tiny segment and has the same period as the myopic path.  Results in this paper are rigorous and inspiring, and the method might help solve the optimal strategy for general games and general algorithms.
\end{abstract}

% keywords can be removed

\section{Introduction}\label{sec1}

This paper studies the problem of the optimal strategy against the learning algorithms in repeated games. The problem can be regarded as a discrete version of the optimal control problem, which is different from the classical optimal control problem. This problem is receiving more attention and becoming more significant from different fields in recent years, especially with the development of artificial intelligence, in which the learning algorithm is the central part.

In the last decade, significant advancements have been made in the field of artificial intelligence (AI), specifically regarding AI algorithms for game playing\cite{silver2016mastering, silver2017mastering, moravvcik2017deepstack, brown2018superhuman, arulkumaran2019alphastar}. Numerous algorithms have been proposed and demonstrated to achieve remarkable performance in various games. Consequently, games involving these algorithms have become increasingly prevalent. This trend is evident in various aspects, such as the widespread participation in AI-powered online games for entertainment, the attempts by malicious individuals to damage face recognition systems\cite{sharif2016accessorize, szegedy2013intriguing, nguyen2015deep}, and the increasing deployment of unmanned systems across diverse industries\cite{finlayson2019adversarial, chen2022unmanned}.  Besides, since algorithms can be considered as approximations of human adaptive behavior, the findings of this study can help decision-making in ubiquitous repeated games in real-world scenarios.

Naturally, this raises fundamental questions: How should we strategize when facing an opponent who utilizes an AI algorithm? How can we design algorithms that yield superior performance even when facing a perfect opponent? Finding answers to these questions is essential for understanding the dynamics of the game between AI and humans. However, the existing literature on this topic lacks sufficient theoretical analysis. Therefore, this paper aims to address the first question within a simplified and ideal scenario.

Algorithms for playing games have been studied for a long time, and compelling achievements have
been made in the past decades. In the literature, games involving algorithms can be roughly classified into two categories: 
games in which all players are algorithms and games in which part of the players are algorithms. 
The analysis for the first category of the game system is called the theory of learning in games and  numerous research has been made\cite{fudenberg1998theory}. In this framework, many classical learning algorithms have been proposed such as fictitious play\cite{robinson1951iterative, miyasawa1961convergence}, no-regret learning algorithm\cite{cesa2006prediction, hart2000simple}. To deal with complex games, recent research tried to design learning algorithms by combining advanced techniques such as neural network\cite{heinrich2015fictitious} and Q-learning\cite{sayin2022fictitious}. 
Then, researchers tried to prove that the time-averaged strategy profile of players can converge to a certain equilibrium\cite{borowski2019learning, papadimitriou2008computing} and different convergence rates have been constructed theoretically\cite{hartline2015no, chen2020hedging}. These studies provided a theoretical basis for the training of game AI. The counterfactual regret minimization (CFR)\cite{zinkevich2007regret} and two-player Texas Hold'em\cite{moravvcik2017deepstack, brown2018superhuman} is a successful instance in this direction. Note that in the research stated above, algorithms are usually assumed to be symmetric or homogeneous. Games between asymmetric algorithms are studied relatively rarely. Bouzy et al.\cite{bouzy2018hedging} proposed the concept of algorithm game and studied the repeated matrix game between different algorithms. By the method of experiment, they got the payoff matrix of the algorithm game and compared the performance of the algorithms especially their performance against the perfect opponent. This is exactly the problem we will study in this paper and our analysis will be theoretical and rigorous.

To find the perfect opponent (or equally the optimal strategy) against these algorithms, one natural idea at the first glance is to take Nash equilibrium strategy at every stage. However, finding Nash equilibrium generally is very hard(actually the problem is PPAD-complete\cite{daskalakis2009complexity, chen2009settling}). Moreover, playing Nash equilibrium at every stage may miss potential benefits since opponents often have hidden behavior patterns in repeated games. This makes it desirable and profitable to deviate from equilibrium and attempt to exploit opponents.

{
There has been some research on the problem of opponent exploitation both experimentally and theoretically, especially in the last decade. Sandholm\cite{sandholm2007perspectives} described a method of tricking the opponent by adopting specific strategies in a few rounds and exploiting him in turn, named the get-taught-and-exploited problem.
Ganzfried and Sandholm\cite{ganzfried2011game} considered how to model and exploit an opponent in large extensive-form games of imperfect information.} Ganzfried et al.\cite{ganzfried2015safe} analyzed how to safely play a finitely-repeated two-player zero-sum game regardless of the strategy used by the opponent. {Ganzfried and Sun\cite{ganzfried2018bayesian} proposed an opponent exploitation method based on the bayesian setting and Dirichlet distribution in imperfect-information games.} Tang et al.\cite{tang2020enhanced} utilized recurrent neural networks to approximate the opponent and the look-ahead strategy to optimize his utility. {These works are basically based on experiments but the experiment cannot reveal the key parameters that affect the optimal control and the internal pattern of the system. Concerning theoretical work,} 
Braverman\cite{braverman2018selling} et al. gave the bounds for the seller’s revenue in an auction where the buyer runs a no-regret learning algorithm; Deng et al.\cite{deng2019strategizing} proved that against a no-regret learner, the player can always guarantee a utility of at least what he would get in a Stackelberg equilibrium of the game under some mild assumptions; Mu and Guo\cite{mu2012towards} considered the problem of the optimal strategy against an opponent with a finite-memory strategy; Dong and Mu\cite{dong2022optimal} studied the optimal strategy against fictitious play in infinitely repeated games. Basically, research in this direction is not sufficient especially when the algorithm is advanced and adaptive.

In this paper, we consider the question that in a finitely repeated two-player zero-sum game, if one player uses the Hedge learning algorithm, how should the other player play in order to maximally exploit it? The Hedge algorithm, also known as multiplicative weight updates\cite{cesa2006prediction}, was generalized by Freund and Schapire\cite{freund1997decision} from the celebrated weighted majority algorithm, which was proposed by Littlestone and Warmuth\cite{littlestone1994weighted}. 
Because of the wide applications and good performance of the Hedge algorithm, exploiting it seems challenging and interesting and we choose it in our study.
First, the Hedge algorithm is ubiquitous in learning theory, game theory and many other areas\cite{plotkin1995fast, kleinberg2009multiplicative, arora2012multiplicative}. Various famous algorithms, such as the AdaBoost algorithm\cite{freund1997decision} in machine learning and the method of mirror descent\cite{wright1999numerical} for convex optimization, employ the core idea of the Hedge algorithm, i.e., the multiplicative weights method. Second, the Hedge algorithm and its variants are widely used in equilibrium computation via self-play in repeated games since the time-averaged strategies of the players converge to coarse correlated equilibrium\cite{hart2000simple} if all players employ no-regret learning algorithms in a general k-person game. In recent years, the Hedge algorithm has been attracting the interest of researchers and some near-optimal learning algorithms have been proposed based on the idea of the Hedge algorithm, such as \cite{daskalakis2011near, piliouras2022beyond}. This is because of the good performance of the Hedge algorithm, i.e., the no-regret nature in online learning with either one sole agent or two symmetric agents in a game\cite{cesa2006prediction}. Thus, the Hedge learning algorithm performs well even faced with strong opponents. Such good adaptation and low exploitability make it more challenging to exploit the algorithm compared with other simpler ones.

The main challenge lies in the coupled dynamics of the game system, i.e., given that one player uses the Hedge algorithm, the strategy of the other player at every stage not only affects the instantaneous return but also affects the future strategies of the Hedge algorithm, thereby affecting his own future returns. To solve the problem, we need to measure the impact of the current strategy on future returns and make a balance between immediate and long-term returns, which is difficult and complex.

Now, we give a brief statement of the problem and results. In this paper, we will consider a $2 \times 2$ zero-sum game repeated for T stages. Suppose one player (player X) employs the Hedge learning algorithm to update his stage strategy. We study the optimal play for the other player (player Y), i.e., how should player Y choose his action sequence(i.e. optimal strategy) in order to obtain the maximal {cumulative expected payoff}?

From the updating formula of the Hedge algorithm, we find that the stage strategy of player X is completely determined by the system {state}. We deduce the updating formula of the state and then construct {the State Transition Triangle Graph}(STTG). Then, the problem of the optimal play is equivalent to the Longest Path Problem on STTG.

First, we analyze the evolution of the game system when player Y adopts {the myopic best response}. We prove that the strategy sequence of player X and the action sequence of player Y enter a cycle after $o(T)$ stages.
Then, we find and prove the recursive rule between the solutions of the Bellman Optimality Equation corresponding to the states at adjacent moments. Based on these results, we further prove that the optimal action sequence of player Y is periodic over the entire repeated game with a small number of stages truncated at the end. The optimal action sequence is different from the myopic action sequence. According to these theoretical results, we provide an efficient and economical method to acquire the optimal play against the Hedge algorithm. 

\textbf{Paper Organization}: we present the problem formulation and assumptions in Section \ref{sec2}; we introduce the state and the State Transition Triangle Graph in Section \ref{sec3}; we investigate the game system when player Y takes myopic best response in Section \ref{sec4}; we study the property of the optimal play and give a simple procedure to obtain it in Section \ref{sec5}; we discuss the impact of the assumptions on the results in Section \ref{sec6}; we conclude this paper and present future work in Section \ref{sec7}.

\section{Problem Formulation}\label{sec2}
Consider a $2\times 2$ zero-sum game. To be specific, there are two players in the game, called players X and Y. Player X has two actions: $U$ and $D$, while player Y has two actions: $L$ and $R$. Given any pure action profile, players X and Y obtain their individual loss/payoff, which can be represented by the following matrix,
\begin{equation}\label{LossMatrix}
 A = \begin{pmatrix}
a_{11} & a_{12}\\
a_{21} & a_{22}
\end{pmatrix}   .
\end{equation}
In the matrix, $a_{11}$ is the loss of player X and the payoff of player Y given the action profile $(U,L)$. For $(U,R),(D,L),(D,R)$, the loss of player X is $a_{12},a_{21}$ and $a_{22}$ respectively. If players X and Y take mixed strategy $x=(x_1,x_2)^T$ and $y=(y_1,y_2)^T$, where $x_1+x_2=1, x_1\geq 0,x_2\geq 0$ and $y_1+y_2=1,y_1\geq 0,y_2\geq 0$, then the expected loss of player X (equally the expected payoff of player Y) is $x^T A y$.

Now let the game repeat $T$ times. Denote the mixed strategy of players X and Y at time $t$ by $x_t=(x_{1,t},x_{2,t})^T$ and $y_t=(y_{1,t},y_{2,t})^T$. 
We assume that player X updates its stage strategy according to the Hedge algorithm. 
Specifically, the initial strategy of player X is set to be $x_1 = (1/2, 1/2)$, and its stage strategy at time $t\geq 2$ is determined by 
\begin{equation}\label{Hedge_initial_formula}
    x_{i,t} = \frac{ \exp(\eta R_{i,t})}{\sum_{j=1}^2\exp(\eta R_{j,t})},\ \ \ i=1,2.
    % {\sum_{j=1}^2 e^{\eta R_{j,t}}}
\end{equation}
The notation $R_{i,t}$ is called \textit{the regret}, defined by
\begin{equation}\label{regret_definition}
R_{i,t} = \sum\limits_{\tau=1}^{t-1} (x_\tau^T A y_\tau-e_i^T A y_\tau).
% = (\hat{L}_T- \min\limits_{i\in\mathcal{E}} { L_{i,T}})
\end{equation}
where $e_1 = (1, 0)^T$ and $e_2 = (0, 1)^T$.
By simple computation, the updating rule \eqref{Hedge_initial_formula} equals
\begin{equation}\label{HedgeStrategyFormula}
    x_{i,t} = \frac{ \exp({-\eta \sum_{\tau=1}^{t-1}e_i^T A y_\tau})}{\sum_{j=1}^2\exp({-\eta \sum_{\tau=1}^{t-1}e_j^T A y_\tau})},\ \ i=1,2.
\end{equation}
In the above formulas, {$\eta$ is called the learning rate} and it is a positive parameter satisfying $\eta = O(\sqrt{\frac1T})$. {In the literature, it is often taken as $\eta = \sqrt{\frac{8 \ln 2}{T}}$ to optimize the regret bound\cite{cesa2006prediction}.}

Now the evolution of the game system depends on the choice of player Y. By choosing different action sequences $y_1,y_2,\cdots,y_T$, player Y can induce different strategy sequences $x_1,x_2,\cdots, x_{T}$ of player X, which results in different \textit{instantaneous expected payoff} (IEP) $x_t^T A y_t$ and \textit{cumulative expected payoff} (CEP) $ \sum_{t=1}^T x_t^T A y_t$ for himself. Then the basic questions arise naturally: 

% i. in order to get the maximal CEP, what action sequence $y_1,y_2,\cdots,y_T$ should player Y choose to get the maximal CEP?

(i). what action sequence should player Y choose to get the maximal CEP?

% ii. under the optimal choice of player Y, how will the system behave?

(ii).  how will the system behave under the optimal choice of player Y?

(iii). can we design a procedure to find the optimal action sequence (also called the optimal strategy) of player Y?

% These are the problems we consider in this paper. For problem i, since the repeated game horizon is finite, we could obtain the optimal action sequence through the Bellman Optimal Equation. However, since solving the Bellman Optimality Equation requires traversing all of the alternative actions of player Y and the reachable states, the cost of this approach explodes as T increases. Therefore, it is necessary and helpful to study the evolution characteristics of the system and further give an economical way to obtain the optimal action sequence, i.e., problems ii and iii.

{
Before presenting the main results, we give some assumptions and remarks on them.
}

{
First, we assume that the elements of the loss matrix are rational, i.e., $a_{ij}\in\mathbb{Q},\ \forall\ i,j=1,2$. This assumption is realistic. Below, in Section \ref{sec3}, Section \ref{sec4}, Section \ref{sec5}, we assume that the elements of the loss matrix are integers for the simplicity of analysis. At the end of Section \ref{sec5}, we will extend the main results to the games with rational elements. Then, in subsection \ref{subsection_irrational_numbers}, we will discuss the situation when some elements are irrational in the loss matrix. }

% In reality, losses are usually integers or rational. When matrix elements are irrational, we can approximate irrational numbers with rational numbers. 

Second, we assume that there is no dominant strategy for player X and player Y.
Let
\begin{align}
    \delta_1=a_{11}-a_{12},\quad  &\delta_2 = a_{21}-a_{22}, \nonumber \\
    \Delta_1 = a_{11}-a_{21},\quad  &\Delta_2 = a_{12}-a_{22}.\nonumber
\end{align}
This assumption means that $\Delta_1\Delta_2<0$ and $\delta_1\delta_2<0$. Without loss of generality, we assume that $\Delta_1 > 0, \Delta_2 < 0, \delta_1>0$ and $\delta_2<0$. For simplicity, we assume that $\vert \Delta_1 \vert < \vert \Delta_2 \vert$, which is not essential.
In subsection \ref{subsection_player_x_has_dominant_strategy}, we study the case when player X has a dominant strategy. When player Y has a dominant strategy and player X does not have a dominant strategy, the problem of the optimal action sequence is much more complex, which is discussed in subsection \ref{subsection_player_Y_has_dominant_strategy}.

Third, we assume that player Y knows that player X uses the Hedge algorithm to update its stage strategy. This seems ideal and unrealistic. However, it is a necessary starting point to study the algorithm game. When player Y does not know the specific algorithm adopted by player X, he can try to identify it and then use the results in this paper. Since the advanced algorithms are often combinations of some simple and classical algorithms such as the Hedge algorithm, results in this paper can provide some basic insights in real plays against these algorithms.

Finally, we assume that player Y adopts pure strategies. Based on STTG, which will be introduced in Section \ref{sec3}, each pure strategy sequence of player Y corresponds to a path on STTG. Hence, each mixed strategy sequence of player Y corresponds to a probability distribution on all paths. Since the game is repeated for finite times, there must exist a path with the longest length, which is just the optimal strategy sequence of player Y and is a pure strategy sequence.

% The maximal exploitation of the Hedge algorithm is essentially an optimal control problem. Since the updating rule of player X is determined in advance, the evolution of the system fully depends on the action sequence of player Y. To get the maximal cumulative returns, player Y needs to determine his strategy sequence to control the behavior of player X before the repeated game starts. If player Y takes mixed strategies, his actual action at each stage depends on random sampling. However, the stage strategy of player X depends on the actual action sequence of player Y rather than the mixed strategy sequence of player Y. Therefore, when player Y adopts mixed strategies, the evolution of the system is not completely controlled by player Y, and thus his cumulative returns are uncertain.

\section{Method: State Transition Triangle Graph}\label{sec3}

In this section, we build the State Transition Triangle Graph based on the updating formula of the Hedge algorithm and IEP of player Y. The problem of the optimal action sequence is proved equivalent to the Longest Path Problem on STTG.

In the updating rule \eqref{HedgeStrategyFormula} of the Hedge algorithm, $x_{1,t}$ and $x_{2,t}$ have the same denominator. Then we have 
\begin{align}
    \frac{x_{2,t}}{x_{1,t}} &= \frac{\exp (-\eta \sum_{\tau = 1}^{t-1}(a_{21},a_{22}) y_\tau)}{\exp (-\eta \sum_{\tau = 1}^{t-1}(a_{11},a_{12}) y_\tau)}\nonumber  \\
    &= \exp(-\eta (-\Delta_1,-\Delta_2) \sum_{\tau = 1}^{t-1} y_\tau)\label{X_2divideX_1}
\end{align}
for time $t\geq 2$.

Let $s_1 = 0$ and
\begin{equation}\label{StateDef}
    s_t \triangleq 
    (-\Delta_1,-\Delta_2) \sum_{\tau = 1}^{t-1} y_\tau, \quad t\geq 2.
\end{equation} 
We call $s_t$ \textit{the state} at time $t$.
Since player Y takes pure strategy at each stage, there are $t$ different possible values of the state at time $t$, denoted by $s_{1,t}, s_{2, t}, \cdots, s_{t, t}$, where 
\begin{equation*}
    s_{i,t} = -(t-i)\Delta_1-(i-1)\Delta_2, 1\leq i\leq t.
\end{equation*}

By \eqref{X_2divideX_1} and \eqref{StateDef}, we have 
\begin{equation}\label{RelationBetweenXandS}
  \frac{x_{2,t}}{x_{1,t}}=\exp(-\eta s_t), \quad t\geq 2.
\end{equation}
This yields 
\begin{equation}\label{NewHedge}
    x_t = \left(\frac{1}{e^{-\eta s_t}+1},\frac{e^{-\eta s_t}}{e^{-\eta s_t}+1}\right), \quad t\geq 2
\end{equation}
since $x_{1,t}+x_{2,t}=1$. When $t=1$, $x_1 = (1/2, 1/2)$ and $s_1 = 0$ also satisfy this relation \eqref{NewHedge}.
Therefore, the stage strategy $x_t$ has a one-to-one correspondence with $s_t$ at all time $t\geq 1$.

%  Given that $s_{t-1}=s_{i,t-1}$, formula \eqref{EqStateUpdatingFormula} can further be written of form
%  \begin{equation}\label{StateTransFormula}
%     s_t = \begin{cases}
%     s_{i,t}, & \text{if}\ y_{t-1} = e_L;\\
%     s_{i+1,t}, & \text{if}\ y_{t-1} = e_R.
%     \end{cases}
% \end{equation}

At time $t$, IEP of player Y is $x_t^T A y_t$, which equals
\begin{equation}\label{PayoffFunc}
    \left(\frac{1}{e^{-\eta s_t}+1},\frac{e^{-\eta s_t}}{e^{-\eta s_t}+1}\right)A y_t \triangleq r(s_t,y_t),
\end{equation}
where $r(s, y)$ is called the payoff function. Note that IEP of player Y at time $t$ only depends on state $s_t$ and action $y_t$.

By \eqref{StateDef}, we have the iteration formula 
\begin{equation}\label{EqStateUpdatingFormula}
    s_{t+1} = s_{t} + (-\Delta_1,-\Delta_2)y_{t} \triangleq h(s_t, y_t)
\end{equation}
and call $h(s, y)$ the state transition function. Hence, state $s_{t+1}$ also only depends on state $s_t$ and action $y_t$.

\begin{figure}[htbp]
    \centering
    \includegraphics[width=0.90\textwidth]{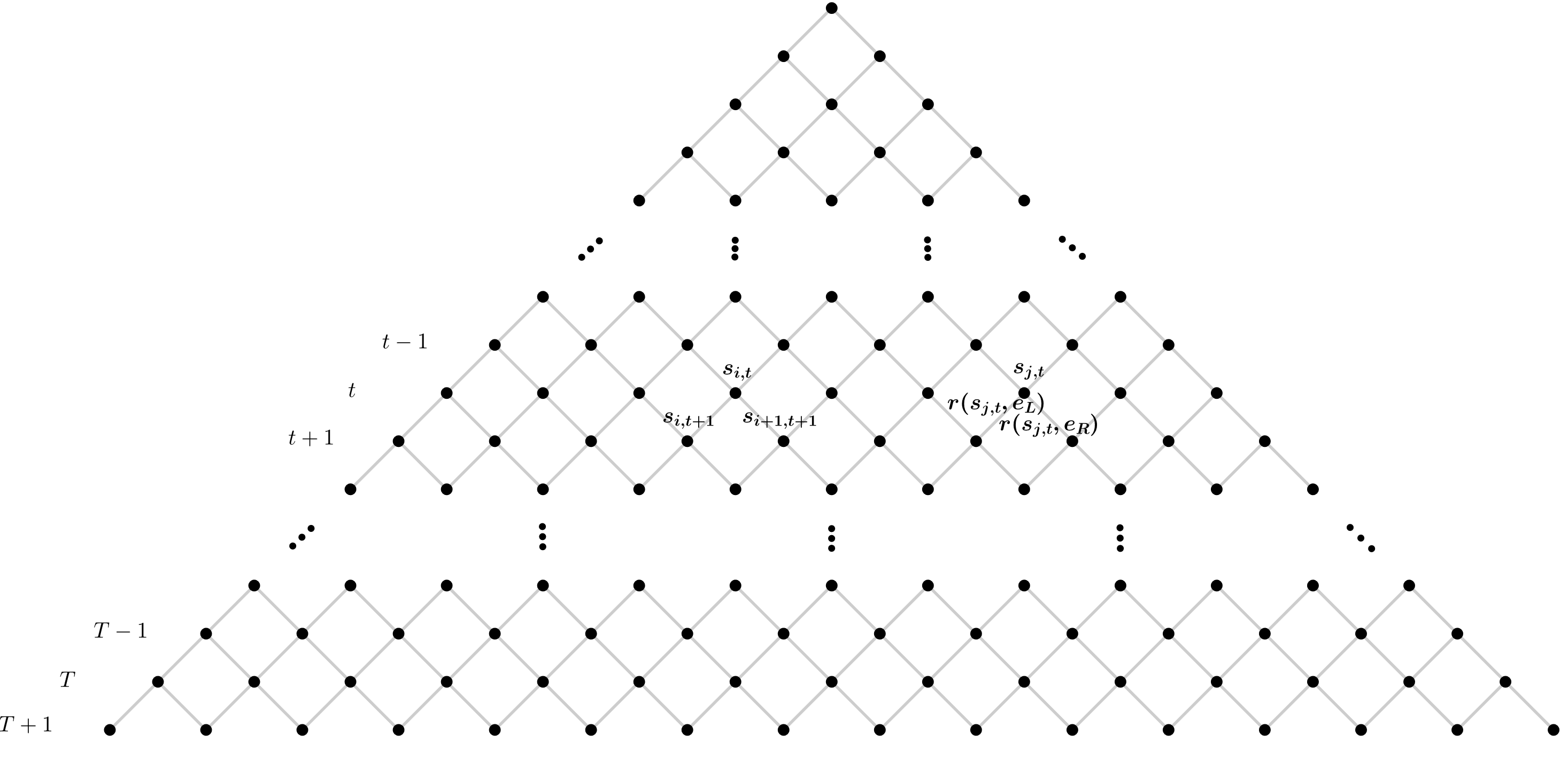}
    \caption{The State Transition Triangle Graph.  }
    \label{fig: AuxiliaryGraphTotal}
\end{figure}

Now, we can build \textit{the State Transition Triangle Graph} (STTG) as shown in Figure \ref{fig: AuxiliaryGraphTotal}, which illustrates the state transition process and IEP of player Y. 

In STTG, each black point, called node, represents a state and the nodes at the $t$-th row represent all possible states at time $t$: $s_{1,t},s_{2,t},\cdots,s_{t,t}$. Each gray edge connecting two nodes represents an action of player Y.
From state $s_{i,t}$, player Y can choose action $e_L$(the left branch), leading to state $s_{i,t+1}$, or chooses action $e_R$(the right branch), leading to state $s_{i+1,t+1}$. Therefore, from state $s_{i, t}$, we have
 \begin{equation}\label{StateTransFormula}
    s_{t+1} = \begin{cases}
    s_{i,t+1}, & \text{if}\ y_{t} = e_L;\\
    s_{i+1,t+1}, & \text{if}\ y_{t} = e_R.
    \end{cases}
\end{equation}
{For each action $y_t\in \{e_L,e_R\}$, IEP $r(s_{i, t}, y_t)$ of player Y can be viewed as a length and is assigned to the according edge. The game is repeated for $T$ times, and the states at time $T+1$ are also drawn for completeness. 
}

Now, we can see that, for player Y, choosing different action sequences is equivalent to choosing different paths on STTG from row $1$ to $T+1$. The CEP of each action sequence of player Y equals the length of the according path. Hence, the problem of the optimal action sequence against the Hedge algorithm equals the Longest Path Problem on STTG. From the STTG (Figure \ref{fig: AuxiliaryGraphTotal}), the feasible paths for player Y increases exponentially with the time horizon T, making it hard to find the optimal one.

\section{Main Results I: Periodicity of the Hedge-myopic System}\label{sec4}

In this section, we study the evolution and performance of the game system, in which player Y takes the myopic best response (MBR) to the stage strategy of player X.

Since player Y knows that player X employs the Hedge algorithm to play the repeated game, player Y can accurately predict the strategy of player X at each stage. It is natural for player Y to take MBR, i.e., 
\begin{equation}\label{myopic_best_response_using_xt}
y_t  =\arg \max\limits_{y} x_t^T A y.
\end{equation}
When player Y adopts this decision rule, the game system is entirely determined. We call the game system \textit{the Hedge-myopic system}. 

Since there is a one-to-one correspondence between $x_t$ and $s_t$, the above equation \eqref{myopic_best_response_using_xt} can be converted to
\begin{equation}\label{myopic_best_response_using_st}
y_t =\arg \max\limits_{y} r(s_t,y) \triangleq \operatorname{MBR}(s_t),
\end{equation}
where $r(s, y)$ is the payoff function defined by \eqref{PayoffFunc}.

Since player Y is assumed to take pure strategy at each stage, to obtain $y_t$, we just need to compare the payoff of action $L$ and $R$ at $s_t$, i.e., $r(s_t, e_L)$ and $r(s_t, e_R)$. By calculating, 
\begin{equation*}
    r(s_t,e_R)-r(s_t,e_L)=-\frac{\delta_1+\delta_2e^{-\eta s_t}}{1+e^{-\eta s_t}}.
\end{equation*}
Since $\delta_1>0$ and $\delta_2<0$, {function $r(s,e_R)-r(s,e_L)$ is monotonically decreasing with respect to $s$ and has a unique zero point}, denoted by
\begin{equation}\label{sast}
    {s^{\ast} = -\frac{1}{\eta}\ln \left(-\frac{\delta_1}{\delta_2}\right).}
\end{equation}
Since $\eta = O(\sqrt{\frac1T})$, this yields $s^\ast = O(\sqrt{T})$. Then, we have $r(s,e_R)\leq r(s,e_L)$ when $s\geq s^\ast$, and $r(s,e_R)>r(s,e_L)$ when $s<s^\ast$.

Now, the function $\operatorname{MBR}$ has a simpler form,
\begin{equation}\label{MBR}
\operatorname{MBR}(s) = 
\begin{cases}
e_R, &\quad \text{if}\ s<s^\ast;\\
e_L, &\quad \text{if}\ s\geq s^\ast.
\end{cases}
\end{equation}

In the Hedge-myopic system, we denote the action sequence of player Y by $y_1^M,y_2^M\cdots,y_{T}^M$, the strategy sequence of player X by $x_1^M,x_2^M,\cdots,x_{T+1}^M$ and the state sequence by $s_1^M,s_2^M,\cdots,s_{T+1}^M$. 
In addition, sequence $s_1^M,s_2^M,\cdots,s_{T+1}^M$ is called \textit{the myopic sequence}, and sequence $s_1^M,y_1^M,s_2^M,y_2^M,\cdots,s_T^M,y_T^M,s_{T+1}^M$ is called \textit{the myopic path}. The evolution process of the Hedge-myopic system is illustrated by Figure \ref{fig:RelationBetweenSequences} below.

\begin{figure}[htbp]
    \centering
    \includegraphics[width=0.90\textwidth]{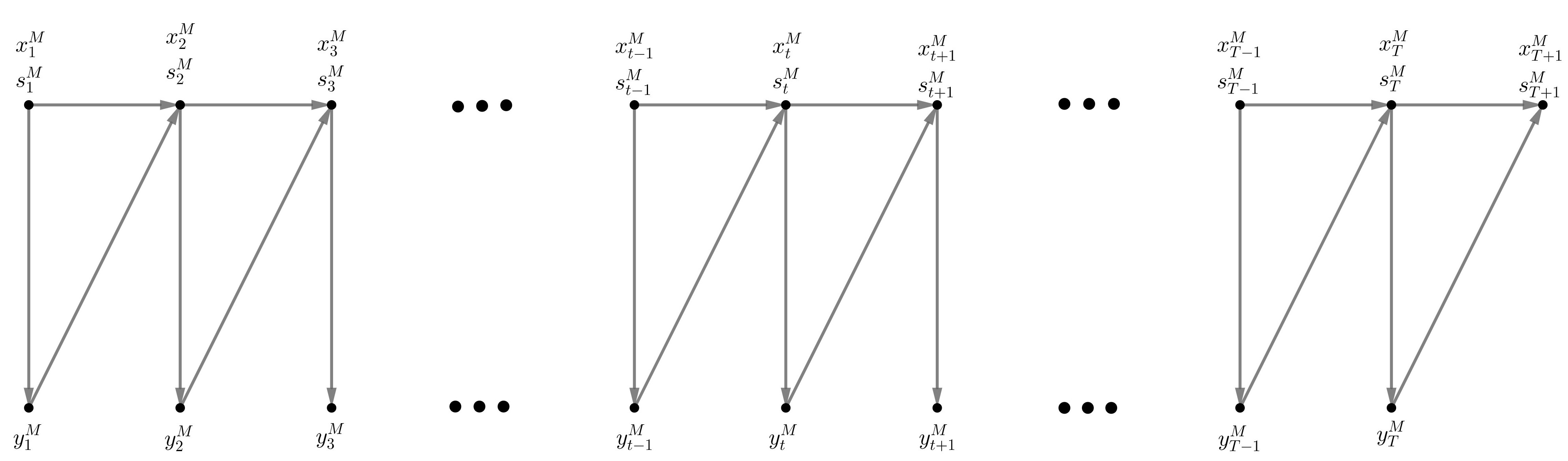}
    \caption{The evolution of the system: Hedge VS. MBR.}
    \label{fig:RelationBetweenSequences}
\end{figure}

From \eqref{MBR}, the states in the Hedge-myopic system satisify
\begin{equation}\label{state_updating_in_myopic_system}
s_{t+1}^M = 
\begin{cases}
s_{t}^M-\Delta_2, &\quad \text{if}\ s_t^M<s^\ast;\\
s_{t}^M-\Delta_1, &\quad \text{if}\ s_t^M\geq s^\ast.
\end{cases}
\end{equation}

% {
% From the state updating formula \eqref{EqStateUpdatingFormula}, we know that $s_{t+1}^M$ is determined by $s_t^M$ and $y_t^M$. By \eqref{MBR}, $y_t^M=\operatorname{MBR}(s_t^M)$ is determined by $s_t^M$. Hence, we have  $s_{t+1}^M$ only depends on $s_t^M$.} The evolution process of the Hedge-myopic system is illustrated by Figure \ref{fig:RelationBetweenSequences}.

By \eqref{state_updating_in_myopic_system}, no matter whether $s_t$ is bigger or smaller than $s^\ast$, it has a trend towards $s^\ast$. This means that if some $s_t^M\geq s^\ast$, the next state $s_{t+1}^M$ will decrease until the state becomes smaller than $s^\ast$; if some $s_t^M<s^\ast$, the next state $s_{t+1}^M$ will increase until the state becomes bigger than $s^\ast$. Hence, the states $s_t^M$ will oscillate around $s^\ast$ and will never move far away from $s^\ast$.

The following theorem reveals the further property of the state sequence $s_t^M$, i.e., it is periodic after a few starting stages, which is our first main result.
\begin{theorem}\label{ThmMyopic}
In the Hedge-myopic system, after $T_s = o(T)$ stages, the action sequence $y_{T_s}^M, y_{T_s+1}^M, \cdots, y_{T}^M$ of player Y, the strategy sequence $x_{T_s}^M, x_{T_s+1}^M, \cdots, x_{T}^M$ of player X and the state sequence $s_{T_s}^M, s_{T_s+1}^M, \cdots, s_{T}^M$ are periodic. Moreover, they have the same least positive period 
$$T^\ast = \frac{m(\vert\Delta_1\vert,\vert\Delta _2\vert)}{\vert\Delta _1\vert}+\frac{m(\vert\Delta _1\vert,\vert\Delta _2\vert)}{\vert\Delta _2\vert},$$
where $m(\vert\Delta _1\vert,\vert\Delta _2\vert)$ is the least common multiple of $\vert\Delta _1\vert$ and $\vert\Delta _2\vert$. 
\end{theorem}

Since $y_t^M$ is the myopic best response to $x_t^M$, each $x_t^M$ determines a unique $y_t^M$. So if $x_{T_s}^M, x_{T_s+1}^M, \cdots, x_{T}^M$ is periodic, then $y_{T_s}^M, y_{T_s+1}^M, \cdots, y_{T}^M$ is also periodic, and they have the same least positive period. On the other hand, $x_t^M$ corresponds to $s_t^M$ one to one, so in order to prove the periodicity of $y_{T_s}^M, y_{T_s+1}^M, \cdots, y_{T}^M$ and $x_{T_s}^M, x_{T_s+1}^M, \cdots, x_{T}^M$, we only need to prove the periodicity of $s_{T_s}^M, s_{T_s+1}^M, \cdots, s_{T}^M$.

We will prove this in several steps. Since $\vert \Delta_2 \vert \geq \vert \Delta_1 \vert$, we find that if the state crosses $s^\ast$ at some point (Lemma \ref{ChangeSign}), there will not exist two consecutive states both less than $s^\ast$. We will prove the periodicity of the subsequence consisting of such states less than $s^\ast$ (Lemma \ref{Lemma4.3}), by studying a new sequence $\{\gamma_i\}_{i=1}^\infty$ generated from the updating rule of the subsequence (Lemma \ref{LemmaNumberTheory}). The logical chain between results in this section is described by Figure \ref{fig:logic_chain_between_lemmas} below.

\begin{figure}[htbp]
    \centering
    \includegraphics[width=0.90\textwidth]{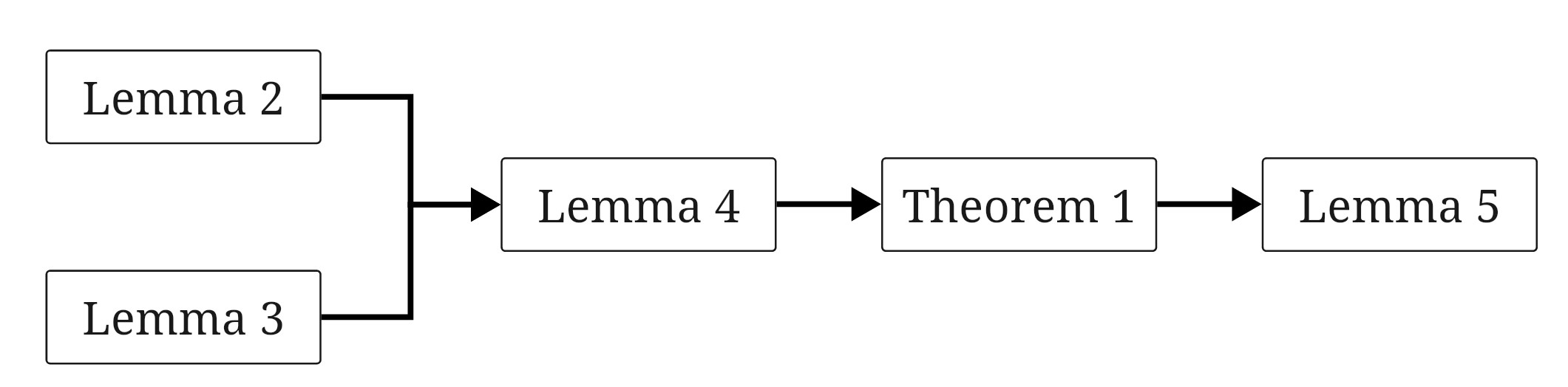}
    \caption{{The logical chain between results in Section \ref{sec4}.}}
    \label{fig:logic_chain_between_lemmas}
\end{figure}

Given any $p,q\in \mathbb{N}$ and $0<p<q$. Given any $\lambda, \gamma \in \mathbb{R}$. Then, there exists a unique integer $k$ such that $\gamma + k p\in [\lambda-p,\lambda)$. Construct the sequence $\{\gamma_i\}_{i=1}^\infty$ by
\begin{equation}\label{GammaSequenceDef}
    \gamma_{i+1} = \gamma_i + q + k_i\cdot p,\quad \forall\ i\geq 1
\end{equation}
where $k_i$ is the unique integer such that $\gamma_{i+1}\in [\lambda-p,\lambda)$. Take $\gamma_1\in [\lambda-p,\lambda)$, then $\gamma_i\in [\lambda-p,\lambda)$ for all $i\geq 1$.

In our proofs, the above notations $p, q, \lambda$ represent $\vert \Delta_1\vert, \vert \Delta_2\vert$ and $s^\ast$ respectively. For the sequence $\{\gamma_i\}_{i=1}^\infty$, we have the following lemma.

\begin{lemma}\label{LemmaNumberTheory}
The sequence $\{\gamma_i\}_{i=1}^\infty$ is periodic, and its least positive period is $\frac{m(p,q)}{q}$.
\end{lemma}

\begin{proof}
By \eqref{GammaSequenceDef}, we have
\begin{align*}
    \gamma_{i+i_0} & = \gamma_{i+i_0-1}+ q+  k_{i+i_0-1}\cdot p\\
     & = \gamma_{i+i_0-2}+ 2q+ k_{i+i_0-2}\cdot p + k_{i+i_0-1}\cdot p\\
     & =\quad  \cdots \\
    & = \gamma_{i}+ i_0\cdot q+\sum\limits_{\tau=i}^{i+i_0-1} k_\tau\cdot p.
\end{align*}
Since $\gamma_i\in [\lambda-p,\lambda)$ for all $i\geq 1$, we obtain that $\gamma_{i+i_0}=\gamma_{i}$ if and only if 
$$i_0\cdot q+\sum_{\tau=i}^{i+i_0-1} k_\tau\cdot p\equiv 0\mod p,$$ 
which is equivalent to $$ i_0\cdot q\equiv 0\mod p.$$ Obviously,  $i_0=\frac{m(p,q)}{q}$ is the least positive integer that satisfies this equation. 
\end{proof}

For the state sequence $s_1^M,s_2^M,\cdots,s_{T+1}^M$, we claim that it will cross $s^\ast$ at some point.

\begin{lemma}\label{ChangeSign}
In the state sequence $s_1^M,s_2^M,\cdots,s_{T+1}^M$, there exists some time $T_c=o(T)$ such that $s_{T_c}^M<s^\ast $ and $s_{T_c+1}^M\geq s^\ast$.
\end{lemma}

\begin{proof}
According to the sign of $s^\ast$, we consider two different cases below. 

\textbf{{Case 1}}: $s^\ast>0$. Since $s_1^M=0$, $s_1^M<s^\ast$. By \eqref{MBR}, $y_1^M=e_R$ and then $$s_2^M=s_1^M+\vert\Delta_2\vert=\vert\Delta_2\vert>s_1^M.$$
If $s_2^M\geq s^\ast$, then we find $T_c$; otherwise, by \eqref{MBR}, $y_2^M=e_R$ and then $s_3^M=s_2^M+\vert\Delta_2\vert=2\vert\Delta_2\vert$. And it goes like this. 

If $s_t^M<s^\ast$ for all $t, 1\leq t\leq T+1$, then when $T$ is sufficiently large, we have $$s_{T+1}^M=T\vert\Delta_2\vert>s^\ast$$ because $s^\ast = o(T)$. Contradiction. Thus, there must exist time $T_c$ such that $s_{T_c}^M<s^\ast,s_{T_c+1}^M\geq s^\ast$. 

\textbf{{Case} 2} : $s^\ast\leq 0$. Since $s_1^M=0$, $s_1^M\geq s^\ast$. By \eqref{MBR}, $y_1^M=e_L$ and then $$s_2^M=s_1^M-\vert\Delta_1\vert=-\vert\Delta_1\vert.$$  
If $s_2^M< s^\ast$, then $y_2^M=e_R$, leading to $$s_3^M=s_2^M+\vert\Delta_2\vert=-\vert\Delta_1\vert+\vert\Delta_2\vert>0\geq s^\ast$$ because $\vert\Delta_1\vert<\vert\Delta_2\vert$. In this case, we find $T_c=2$. If $s_2^M\geq s^\ast$, by \eqref{MBR}, $y_2^M=e_L$ and then $$s_3^M=s_2^M-\vert\Delta_1\vert=-2\vert\Delta_1\vert.$$ And it goes like this. 

If $s_t^M\geq s^\ast$ for all $ 1\leq t\leq T+1$, then when $T$ is sufficiently large, we have $$s_{T+1}^M=-T\vert\Delta_1\vert<s^\ast$$ because $s^\ast = o(T)$. Contradiction. Thus, there must exist time $t_c$ such that $s_{t_c}^M\geq s^\ast,s_{t_c+1}^M< s^\ast$. By \eqref{MBR}, $y_{t_c+1}^M = e_R$ and then $s_{t_c+2}^M=s_{t_c}^M-\vert\Delta_1\vert+\vert\Delta_2\vert>s_{t_c}^M\geq s^\ast$. In this case, we find $T_c=t_c+1$.

From the above discussion of $T_c$, we can see that the least $T_c$ satisfies
\begin{equation*}
    T_c\leq \max \{\lceil \frac{\vert s^\ast\vert}{\vert \Delta_1 \vert}\rceil,\lceil \frac{\vert s^\ast\vert}{\vert \Delta_2 \vert}\rceil \} +2 =\max \{ \lceil \frac{\vert s^\ast\vert}{\vert \Delta_1 \vert}\rceil\} +2.
\end{equation*}
Thus, $T_c = o(T)$ because $s^\ast = o(T)$.
\end{proof}

This lemma says that the state will cross $s^\ast$ from the side less than $s^\ast(T)$ to the side greater than $s^\ast$ at some time $T_c$.
Denote the least $T_c$ by $t_0$, then $t_0 = o(T)$. We define a time sequence by
\begin{equation}\label{Ti+1Def}
    t_{i+1} = \min \{t: t_i<t\leq T,s_{t}^M<s^\ast,\ s_{t+1}^M>s^\ast\}.
\end{equation}
Denote $t_{i_{\max}}$ the maximal time $t$ satisfying $t\leq T,$ $s_{t}^M<s^\ast$ and $s_{t+1}^M>s^\ast$.
Then, the state sequence $s_{t_1}^M,s_{t_2}^M,\cdots,s_{t_{i_{\max}}}^M$ consists of all states less than $s^\ast$. 

The following lemma shows that the state subsequence is periodic and its least positive period is $T^\ast$.

\begin{lemma}\label{Lemma4.3}
For the sequence $s_{t_1}^M,s_{t_2}^M,\cdots,s_{t_{i_{\max}}}^M$, we have
\begin{equation*}
    s_{t_i+T^\ast}^M = s_{t_i}^M, \quad \forall\ 1\leq i\leq i_{\max},
\end{equation*}
where $T^\ast = \frac{m(\vert\Delta_1\vert,\vert\Delta_2\vert)}{\vert\Delta_1\vert}+\frac{m(\vert\Delta_1\vert,\vert\Delta_2\vert)}{\vert\Delta_2\vert}$.
\end{lemma}

\begin{proof}
By \eqref{Ti+1Def},  $s_{t_i}^M<s^\ast$ and then by \eqref{MBR}, $\operatorname{MBR}(s_{t_i}^M)=e_R$, leading to $s_{t_i+1}^M = s_{t_i}^M+\vert\Delta_2\vert$. Again by \eqref{Ti+1Def}, $s_{t_i+1}^M> s^\ast$ and by \eqref{MBR}, $\operatorname{MBR}(s_{t_i+1}^M)=e_L$. 

Now, consider the states $s_{t_i+1}^M,s_{t_i+2}^M,\cdots,s_{t_{i+1}-1}^M$ where $0\leq i<i_{\max}$. We prove that all of them are no less than $s^\ast$, i.e., 
\begin{equation}\label{GeqsAst}
    s_{t_i+j}^M\geq s^\ast
\end{equation}
for all $t_i+1\leq t_i+j\leq t_{i+1}-1$.

If not, there exists time $t_i+j$, $t_i+1< t_i+j\leq t_{i+1}-1$ such that $s_{t_i+j-1}^M\geq s^\ast,s_{t_i+j}^M<s^\ast$. Then, by \eqref{MBR}, $\operatorname{MBR}(s_{t_i+j-1}^M)=e_L$, $\operatorname{MBR}(s_{t_i+j}^M)=e_R$ and $s_{t_i+j+1}^M=s_{t_i+j-1}^M-\vert\Delta_1\vert+\vert\Delta_2\vert>s_{t_i+j-1}^M\geq s^\ast$ since $\vert\Delta_1\vert<\vert\Delta_2\vert$. Thus, we have $s_{t_i+j}^M<s^\ast$, $s_{t_i+j+1}^M>s^\ast$ which contradicts with the definition \eqref{Ti+1Def}  of $t_{i+1}$.

Hence, we have $s_{t_i+j}^M\geq s^\ast$ for all $t_i+j$ where $t_i+1\leq t_i+j\leq t_{i+1}-1$. Specially, $s_{t_i-1}^M\geq s^\ast$ for all $1\leq i\leq i_{\max}$.

% This contradicts the definition of of , by \eqref{Ti+1Def}, $t_{i+1}=t_i+j<t_{i+1}$. Contradiction.

On one hand,  by \eqref{MBR}, we have $\operatorname{MBR}(s_t^M) = e_L$ for $t=t_i+1,\cdots,t_{i+1}-1$, leading to
 \begin{equation*}
     \begin{aligned}
     s_{t_{i+1}}^M&=s_{t_i}^M-\Delta_2-(t_{i+1}-t_i)\Delta_1\\
     & = s_{t_i}^M+\vert\Delta_2\vert+(t_{i}-t_{i+1})\vert\Delta_1\vert.
     \end{aligned}
 \end{equation*}

On the other hand, we have $s_{t_i-1}^M\geq s^\ast$, $s_{t_i}^M<s^\ast$ for all $i\geq 1$.
By \eqref{MBR}, $\operatorname{MBR}(s_{t_i-1}^M)=e_L$ and then 
$s_{t_i}^M=s_{t_i-1}^M-\vert\Delta_1\vert\geq s^\ast-\vert\Delta_1\vert$. So we have $$s_{t_i}^M\in [s^\ast-\vert\Delta_1\vert,s^\ast)$$ for all $i\geq 1$.

Recall \eqref{GammaSequenceDef} and Lemma \ref{LemmaNumberTheory}, we know that the sequence $s_{t_1}^M,s_{t_2}^M,\cdots, s_{t_{i_{\max}}}^M$ is periodic and its period is $i^\ast = \frac{m(\vert\Delta_1\vert,\vert\Delta_2\vert)}{\vert\Delta_2\vert}$, i.e.,
\begin{equation}\label{EqTemp1}
    s_{t_{i+i^\ast}}^M=s_{t_i}^M.
\end{equation}

Now, consider the state sequence $s_{t_i}^M,s_{t_{i}+1}^M,\cdots,s_{t_{i+i^\ast}-1}^M$. By \eqref{GeqsAst}, 
we know that in this sequence, $s_t^M<s^\ast$ only for $t=t_i,t_{i+1},\cdots,t_{i+i^\ast-1}$, that is to say, $\operatorname{MBR}(s_t^M)=e_R$ only for $t=t_i,t_{i+1},\cdots,t_{i+i^\ast-1}$. Thus, in the sequence $y_{t_i}^M,y_{t_{i}+1}^M,\cdots,y_{t_{i+i^\ast}-1}^M$, the number of $e_R$, denoted by $\operatorname{num}(e_R)$, is
\begin{equation}\label{EqTemp2}
    \operatorname{num}(e_R) = i^\ast = \frac{m(\vert\Delta_1\vert,\vert\Delta_2\vert)}{\vert\Delta_2\vert}.
\end{equation}

By the state updating formula \eqref{EqStateUpdatingFormula}, we have
\begin{equation}\label{EqTemp3}
        \begin{aligned}
        s_{t_{i+i^\ast}}^M&= s_{t_i}^M+\sum_{\tau=t_i}^{t_{i+i^\ast}-1}(-\Delta_1,-\Delta_2)y_\tau\\
        &= s_{t_i}^M-\operatorname{num}(e_L)\vert\Delta_1\vert+\operatorname{num}(e_R)\vert\Delta_2\vert,
        \end{aligned}
\end{equation}
where $\operatorname{num}(e_L)$ is the number of $e_L$ in the sequence $y_{t_i}^M,y_{t_{i}+1}^M,\cdots,y_{t_{i+i^\ast}-1}^M$.

By \eqref{EqTemp1},\eqref{EqTemp2},\eqref{EqTemp3}, we have
\begin{equation*}
    \operatorname{num}(e_L)=\frac{m(\vert\Delta_1\vert,\vert\Delta_2\vert)}{\vert\Delta_1\vert}.
\end{equation*}
So we have $t_{i+i^\ast}-t_i = \frac{m(\vert\Delta_1\vert,\vert\Delta_2\vert)}{\vert\Delta_2\vert}+\frac{m(\vert\Delta_1\vert,\vert\Delta_2\vert)}{\vert\Delta_1\vert}= T^\ast$ and 
\begin{equation*}
    s_{t_i+T^\ast}^M=s_{t_{i+i^\ast}}^M=s_{t_i}^M.
\end{equation*}
\end{proof}

Now, we have proved that a subsequence of $s_t^M$ is periodic. As we have stated, $s_{t+1}^M$ only depends on $s_t^M$. Combing these leads us to prove Theorem \ref{ThmMyopic} as below. 

\begin{proof}
By Lemma \ref{Lemma4.3}, the state subsequence $s_{t_1}^M,s_{t_2}^M,\cdots,s_{t_{i_{\max}}}^M$ is periodic and its least positive period is $T^\ast = \frac{m(\vert\Delta_1\vert,\vert\Delta _2\vert)}{\vert\Delta _1\vert}+\frac{m(\vert\Delta _1\vert,\vert\Delta _2\vert)}{\vert\Delta _2\vert}$. Hence, take $T_s=t_1$, and then the state sequence $s_{T_s}^M, s_{T_s+1}^M, \cdots, s_{T}^M$ is periodic and its least positive period is $T^\ast = \frac{m(\vert\Delta_1\vert,\vert\Delta _2\vert)}{\vert\Delta _1\vert}+\frac{m(\vert\Delta _1\vert,\vert\Delta _2\vert)}{\vert\Delta _2\vert}.$
\end{proof}

\begin{remark}
    In the Hedge-myopic system, the next state $s_{t+1}^M$ only depends on the current state $s_t^M$. If $\operatorname{MBR}(s_t^M) = e_L$, then $s_{t+1}^M = s_t^M - \Delta_1$. If $\operatorname{MBR}(s_t^M) = e_R$, then $s_{t+1}^M = s_t^M - \Delta_2$. Since $\Delta_1\Delta_2<0$ and are both rational, there exist two integers $p$ and $q$ such that $p \Delta_1 + q\Delta_2 = 0$. This explains why the sequence $s_t^M$ is periodic after a few stages. If $\Delta_1$ or $\Delta_2$ is irrational, then the periodicity may not hold.
\end{remark}

\begin{remark}
    From the proof of Lemma \ref{Lemma4.3}, we can see that there are $\frac{m(\vert\Delta_1\vert,\vert\Delta_2\vert)}{\vert\Delta_1\vert}$ action $L$ and $\frac{m(\vert\Delta_1\vert,\vert\Delta_2\vert)}{\vert\Delta_2\vert}$ action $R$ in a period. Then, the time-averaged strategy of player Y in a period is $\left(\frac{\vert \Delta_2\vert}{\vert \Delta_2\vert+\vert \Delta_1\vert}, \frac{\vert \Delta_1\vert}{\vert \Delta_2\vert+\vert \Delta_1\vert}\right)$, which is exactly the Nash equilibrium strategy. However, player Y can obtain higher average payoff than game value now since at each stage player Y is taking a targeted action against the stage strategy of player X. 
\end{remark}

Next, we give an example to illustrate the above results. 

\begin{figure}[htbp]
    \centering
    \includegraphics[width = 0.9\textwidth]{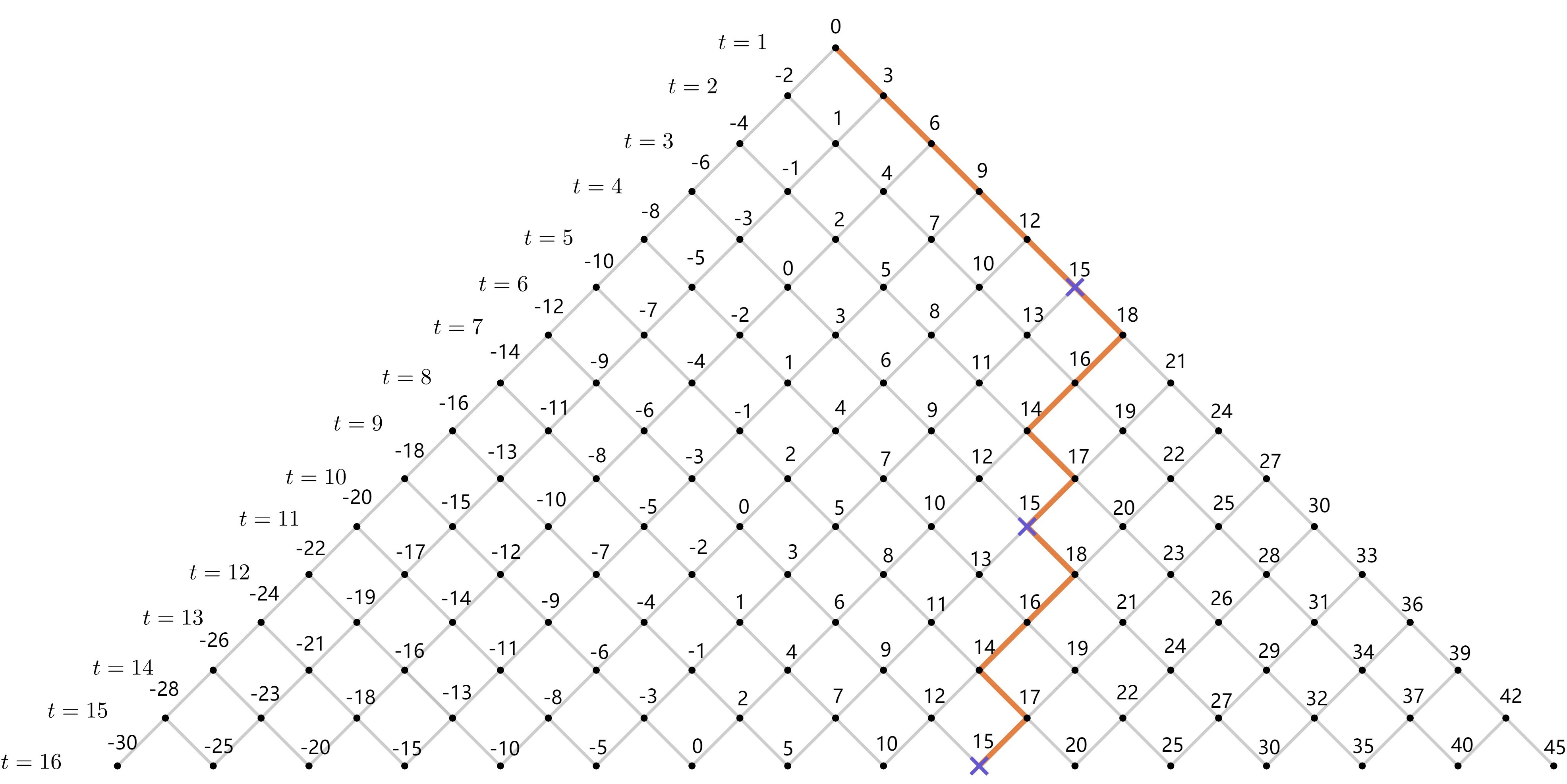}
    \caption{STTG with the myopic path when $\Delta_1 = 2, \Delta_2 = -3, s^\ast = 15.58$. Only the first 16 stages are drawn. The numbers over the nodes represent the values of states $s_t^M$. The orange polyline represents the myopic path.}
    \label{fig:exampleMyopic}
\end{figure}

\begin{example}\label{EX1}
Given that $T=700$ and \begin{equation*}
    A = \begin{pmatrix}
    1 & 0\\
    -1 & 3
    \end{pmatrix}.
\end{equation*}
By simple calculation, 
$\Delta_1 = 2,\Delta_2=-3,\delta_1 = 1,\delta_2 = -4,s^\ast \approx 15.58$. STTG for this example are illustrated by Figure \ref{fig:exampleMyopic}.  In this case, sequence $s_{t_0}^M, s_{t_1}^M,\cdots,s_{t_{i^{\max}}}^M$ is $s_6^M, s_9^M,s_{11}^M,s_{14}^M$ and from time $t=9$, sequence $s_9^M,s_{10}^M,\cdots,s_{16}^M$ (i.e. 14, 17, 15, 18, 16, 14, 17, 15) is periodic and its least positive period is $5$. Then the myopic path can be obtained and is represented by the orange line in Figure \ref{fig:exampleMyopic}.
In each period, the action sequence of player Y is $R, L, L, R, L$ and thus time-averaged strategy of player Y over one period is $(0.6,0.4)$. On the other hand, in one period, payoff sequence obtained by player Y is 0.6646, 0.6119, 0.6702, 0.6390, 0.6250 and thus the time-averaged payoff of player Y is 0.6421, which is higher than the game value 0.6 .
\end{example}

Now, consider a special kind of game where the Nash equilibrium strategy of player X is $(\frac{1}{2},\frac{1}{2})$. In this case, we have $\delta_1  =-\delta_2$ and $s^\ast = -\frac{1}{\eta}\ln(-\frac{\delta_1}{\delta_2}) = 0$. 

For this game, the myopic state sequence is called \textit{the zero sequence} and denoted by $s_1^0,s_2^0,\cdots,s_{T+1}^0$. The according action sequence of player Y is denoted by $y_1^0,y_2^0,\cdots,y_T^0$ and sequence $s_1^0,y_1^0,s_2^0,y_2^0,\cdots,s_T^0,y_T^0,s_{T+1}^0$ is called \textit{the zero path}. Once $\Delta_1$ and $\Delta_2$ is given, the zero sequence and the zero path can be obtained. 

\begin{lemma}\label{LemmaProperties}For the zero sequence, we have

(i) $s_1^0,s_2^0,\cdots,s_{T+1}^0$ is periodic and its least positive period is $T^\ast=\frac{m(\vert\Delta_1\vert,\vert\Delta_2\vert)}{\vert\Delta_2\vert}+\frac{m(\vert\Delta_1\vert,\vert\Delta_2\vert)}{\vert\Delta_1\vert}$.

(ii) If $s_t^0<0$, then $s_{t-1}^0\geq 0$.

(iii) $-\vert\Delta_1\vert\leq s_{t}^0<\vert \Delta_2 \vert$.

(iv) For possible value $s_{i,t}$ of state at time $t$, if $s_{i,t}<0$ and $s_{i+1,t}\geq 0$, then $s_{t+1}^0=s_{i+1,t+1}$.
% \begin{enumerate}
%     \item $s_1^0,s_2^0,\cdots,s_{T+1}^0$ is periodic and its least positive period is $T^\ast=\frac{m(\vert\Delta_1\vert,\vert\Delta_2\vert)}{\vert\Delta_2\vert}+\frac{m(\vert\Delta_1\vert,\vert\Delta_2\vert)}{\vert\Delta_1\vert}$.
%     \item If $s_t^0<0$, then $s_{t-1}^0\geq 0$.
%     \item $-\vert\Delta_1\vert\leq s_{t}^0<\vert \Delta_2 \vert$.
%     \item For possible value $s_{i,t}$ of state at time $t$, if $s_{i,t}<0$ and $s_{i+1,t}\geq 0$, then $s_{t+1}^0=s_{i+1,t+1}$.

% \end{enumerate}
\end{lemma}

\begin{proof}We prove Lemma \ref{LemmaProperties} as follows:

(i) By \eqref{MBR}, $\operatorname{MBR}(s_1)=e_L$. Then, $s_2^0=s_1^0-\vert\Delta_1\vert=-\vert\Delta_1\vert<0$. Again by \eqref{MBR}, $\operatorname{MBR}(s_2^0)=e_R$, leading to $s_3^0=-\vert\Delta_1\vert+\vert\Delta_2\vert>0$ since $\vert\Delta_1\vert<\vert\Delta_2\vert$. 

Take $t_0=2$. Define $t_{i+1} = \min \{t:t_i<t\leq T+1,s_t^0<0,s_{t+1}^0>0\}$. Since $s_2^0=-\vert\Delta_1\vert\in [-\vert\Delta_1\vert,0)$, recall the proof of Theorem \ref{ThmMyopic}, we know that sequence $s_{t_0}^0,s_{t_0+1}^0,\cdots, s_{T+1}^0$ is periodic and its least positive period is $T^\ast = \frac{m(\vert\Delta_1\vert,\vert\Delta_2\vert)}{\vert\Delta_2\vert}+\frac{m(\vert\Delta_1\vert,\vert\Delta_2\vert)}{\vert\Delta_1\vert}$.

Specially, $s_{2+T^\ast}^0 = s_2^0<0$. By \eqref{GeqsAst}, we have $s_{2+T^\ast-1}^0=s_{1+T^\ast}^0\geq 0$. Then, 
by \eqref{MBR}, $\operatorname{MBR}(s_{1+T^\ast}^0)=e_R$ and $s_{2+T^\ast}^0 = s_{1+T^\ast}^0-\vert\Delta_1\vert$. So we have
$s_{1+T^\ast}^0 = s_{2+T^\ast}^0+\vert\Delta_1\vert = s_{2}^0+\vert\Delta_1\vert = s_{1}^0 = 0$. 
This proves that the full sequence $s_1^0,s_2^0,\cdots,s_{T+1}^0$ is periodic.

(ii) By \eqref{GeqsAst}, it obviously holds. 

(iii) 
The minimal $s_t^0$ must be less than $0$, i.e., $s_t^0<0$. Otherwise, by \eqref{MBR}, $\operatorname{MBR}(s_t^0)=e_L$ and then $s_{t+1}^0=s_{t}^0-\vert\Delta_1\vert<s_{t}^0$, which contradicts to the assumption that $s_t^0$ is the minimal state. Then, by \eqref{GeqsAst}, $s_{t-1}^0\geq 0$ and $s_t^0=s_{t-1}^0-\vert \Delta_1\vert \geq -\vert \Delta_1\vert $. 

On the other hand, the maximal $s_t^0$ must satisfy that $s_t^0\geq 0$. Otherwise, $s_t^0<0$ and then by \eqref{GeqsAst}, $s_{t-1}^0\geq 0>s_t^0$. Contradiction. 

And for time $t$ when $s_t^0$ takes the maximum, we have $s_{t-1}<0$. If not, $s_{t-1}^0\geq 0$ and by \eqref{MBR}, we have $s_{t}^0=s_{t-1}^0-\vert \Delta_1\vert <s_{t-1}^0$. Contradiction.

Therefore, we have $s_{t-1}^0<0$ and $s_t^0\geq 0$. Then, by \eqref{MBR}, we have $\operatorname{MBR}(s_{t-1}^0) = e_R$ and then $s_t^0=s_{t-1}^0+\vert \Delta_2\vert <\vert \Delta_2\vert $.

(iv) If $s_{i,t}<0$ and $s_{i+1,t}\geq 0$, $s_{i-1,t}=s_{i,t}-\vert\Delta_1\vert-\vert\Delta_2\vert<-\vert\Delta_1\vert$ and $s_{i+2,t}=s_{i+1,t}+\vert\Delta_1\vert+\vert\Delta_2\vert>\vert\Delta_2\vert$. Since we have proven that $-\vert\Delta_1\vert\leq s_{t}^0<\vert \Delta_2 \vert$, then either $s_t^0=s_{i,t}$ or $s_t^0=s_{i+1,t}$. 

If $s_t^0=s_{i,t}<0$, by \eqref{MBR}, $\operatorname{MBR}(s_{i,t})=e_R$ and then by \eqref{StateTransFormula}, $s_{t+1}^0=s_{i+1,t+1}$. If $s_t^0=s_{i+1,t}\geq 0$, by \eqref{MBR}, $\operatorname{MBR}(s_{i+1,t})=e_L$ and then by \eqref{StateTransFormula}, $s_{t+1}^0=s_{i+1,t+1}$. 
In both cases, $s_{t+1}^0=s_{i+1,t+1}$.
\end{proof}

By Lemma \ref{LemmaProperties}, the zero path, a special myopic path, is periodic from the beginning stage. This is different from the general myopic path, which is periodic after a few stages by Theorem \ref{ThmMyopic}. Thus, the myopic path for a game depends not only on $\Delta_1$ and $\Delta_2$, but also on $\delta_1$ and $\delta_2$, i.e., all the loss elements. To be specific, $\delta_1$ and $\delta_2$ determine $s^\ast$, around which the myopic path oscillate; $\Delta_1$ and $\Delta_2$ determine the period of the myopic path.

The periodicity of the Hedge-myopic system inspires us to investigate whether the nice property still holds for the game system under the optimal action sequence of player Y. {In the next section, we will prove that the answer is basically positive. }

\section{Main Results II: the Optimal Play against Hedge}\label{sec5}

In this section, we will present our main result about the optimal play of player Y. Since the game is repeated $T$ times, theoretically we can use the Bellman Optimality Equation to solve the optimal play.

By STTG in Figure \ref{fig: AuxiliaryGraphTotal}, we can see that from state $s_{i,t}$, there are only two possible states for time $t+1$: $s_{i,t+1}$ and $s_{i+1,t+1}$, which implies that from any state $s_{t}$, not all possible states at time $t^\prime>t$ can be reached. Thus, for $t^\prime >t$, state $s_{t^\prime }$ is said to be \textit{accessible} from state $s_t$ if there is an action sequence $y_t,\cdots,y_{t^\prime-1}$ of player Y such that $s_{t^\prime} = s_{t}+\sum_{\tau=t}^{t^\prime -1}y_\tau.$ We denote such sequence $y_t,\cdots,y_{t^\prime-1}$ by $p(s_t,s_{t^\prime})$. Note that $p(s_t,s_{t^\prime})$ may not be unique, so denote the set of $p(s_t,s_{t^\prime})$ by $\mathcal{P}(s_t,s_{t^\prime })$.

Denote $\mathcal{Z}(s_t)$ to be the set of states at the terminal time $T+1$ which is accessible from state $s_t$.
 
Given $s_t$ and $s_{t^\prime}$ where $s_{t^\prime}$ is accessible from $s_{t}$, define 
\begin{equation*}
  f^\ast(s_t,s_{t^\prime}) =\max\limits_{p(s_t,s_{t^\prime})\in \mathcal{P}(s_t,s_{t^\prime})} \sum_{\tau=t}^{t^\prime-1}r(s_\tau,y_\tau).
\end{equation*}
and $$f^\ast(s_t) = \max\limits_{s_{T+1}\in \mathcal{Z}(s_t)} f^\ast(s_t,s_{T+1}).$$
Then, we have the Bellman Optimality Equation for our problem:
\begin{equation}\label{EqBellmanOptimalEquation}
f^\ast(s_t) = \max\limits_{y\in \{e_L,e_R\}}\{r(s_{t},y)+f^\ast(h(s_t,y))\},\quad 1\leq t\leq T
\end{equation}
where $f^\ast(h(s_T,y))\equiv 0$. 

Define
\begin{equation}\label{yast}
   y^\ast(s_t) = \arg\max\limits_{y\in \{e_L,e_R\}}\{r(s_{t},y)+f^\ast(h(s_t,y))\},
\end{equation}
which is called \textit{the optimal action} of state $s_t$.

Our goal is to find the action sequence of player Y that maximizes $f^\ast(s_1)$. The optimal action sequence is called \textit{the optimal play} and denoted by $y_1^\ast,\ldots,y_{T}^\ast$. We denote the according state sequence by $s_1^\ast,s_2^\ast,\cdots,s_{T+1}^\ast$ and call the sequence $s_1^\ast,y_1^\ast,s_2^\ast,y_2^\ast,\cdots,s_T^\ast,y_T^\ast,s_{T+1}^\ast$ \textit{the optimal path}.

To obtain the optimal path, we need to solve the Bellman Optimality Equation \eqref{EqBellmanOptimalEquation}. However, solving the Bellman Optimality Equation requires traversing all the alternative actions of player Y and all the states on STTG, and thus the cost explodes as T increases. Therefore, it is necessary to study the further internal characteristic of the game system, by which we can answer the problem (ii) and (iii) proposed in Section \ref{sec2} at the same time. Before presenting the analysis in detail, we first give the main result in the theorem below.

\begin{theorem}\label{ThmBest}
There exists time $t^d$ such that the optimal play $y_1^\ast,\cdots,y_{t^d}^\ast$ of player Y is periodic and $T-t^d = o(T)$. Moreover, its least positive period is $T^\ast=\frac{m(\vert\Delta_1\vert,\vert\Delta_2\vert)}{\vert\Delta_2\vert}+\frac{m(\vert\Delta_1\vert,\vert\Delta_2\vert)}{\vert\Delta_1\vert}$.
\end{theorem}

\begin{figure}[htbp]
    \centering
    \includegraphics[width = 0.9\textwidth]{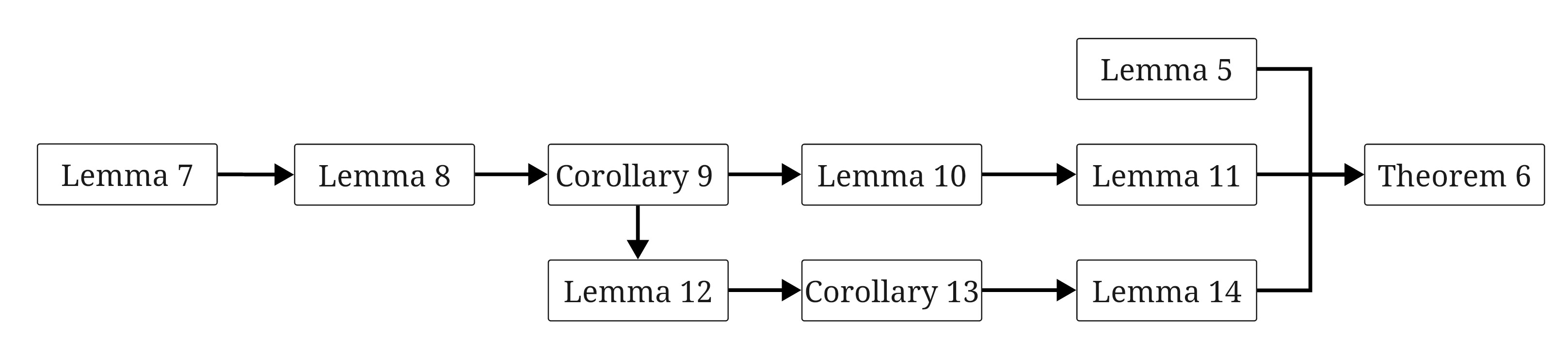}
    \caption{{The logical chain between the results in Section \ref{sec5}.}}
    \label{fig:logic_chain_in_section_5}
\end{figure}

To prove Theorem \ref{ThmBest}, we basically need to solve the Bellman Optimality Equation. Due to the nature of the problem, we can get the recurrence relation of the optimal solution of the Bellman Equation (Lemma \ref{Lemmas0Difference} and Corollary \ref{LemmaRtoRLtoL}) and build a connection between the zero path and the optimal path (Lemma \ref{Lemma1}). Then, we can solve the Bellman Equation in several steps (Lemma \ref{Lemma2}, Lemma \ref{Lemma3}, Corollary \ref{Corollary1} and Lemma \ref{Lemma4}). The logical chain between the results in this section is illustrated in Figure \ref{fig:logic_chain_in_section_5}.

% Applying backward induction, the Bellman Optimality Equation \eqref{EqBellmanOptimalEquation} can be solved. However, the cost of this method explodes when $T$ is too large. Therefore, we wish to explore the specific inherent features of STTG to find a more efficient method to solve the optimal play. 

\subsection{Recurrence Relation of the Optimal Action}

From STTG, we can see that from any state $s$, taking action $e_L$ and then $e_R$ leads to the same state with taking action $e_R$ and then $e_L$. 
However, they may lead to different cumulative expected payoffs of player Y, i.e.,  $r(s,e_R)+r(s-\Delta_2,e_L) \neq  r(s,e_L)+r(s-\Delta_1,e_R)$. 

Define 
\begin{gather*}
    P_{LR}(s) \triangleq r(s,e_L)+r(s-\Delta_1,e_R),\\
    P_{RL}(s) \triangleq r(s,e_R)+r(s-\Delta_2,e_L).
\end{gather*}

Next, we study the difference between $P_{{LR}}(s)$ and $P_{{RL}}(s)$, which turns out to be crucially important in the proofs of the following results. 
By \eqref{PayoffFunc}, we have
\begin{equation*}
\begin{aligned}
    P_{R L}(s)-P_{LR}(s) &= \frac{-\Delta_2(1-e^{\eta \Delta_1})}{(e^{\eta s}+e^{\eta \Delta_1})(1+e^{-\eta s})}+\frac{\Delta_1(1-e^{\eta \Delta_2})}{(e^{\eta s}+e^{\eta \Delta_2})(1+e^{-\eta s})}\\
    & = \frac{-\Delta_2(1-e^{\eta \Delta_1})(e^{\eta s}+e^{\eta \Delta_2})+\Delta_1(1-e^{\eta \Delta_2})(e^{\eta s}+e^{\eta \Delta_1})
    }{(e^{\eta s}+e^{\eta \Delta_1})(e^{\eta s}+e^{\eta \Delta_2})(1+e^{-\eta s})}.
\end{aligned}
\end{equation*}
Let
\begin{align}
         f(s)&\triangleq -\Delta_2(1-e^{\eta \Delta_1})(e^{\eta s}+e^{\eta \Delta_2})+\Delta_1(1-e^{\eta \Delta_2})(e^{\eta s}+e^{\eta \Delta_1}) \label{f(s)Def}  \\ 
        & = \left(\Delta_1(1-e^{\eta\Delta_2})-\Delta_2(1-e^{\eta \Delta_1})\right)e^{\eta s}+(\Delta_1(1-e^{\eta\Delta_2})e^{\eta \Delta_1}  \nonumber\\ 
         &\ \ -\Delta_2(1-e^{\eta \Delta_1})e^{\eta \Delta_2})\nonumber \\
        & \triangleq f_1\cdot e^{\eta s}+f_2\nonumber   
\end{align}
where $f_1$ and $f_2$ are two constants. Obviously, the sign of $f(s)$ is the same as the sign of $P_{R L}(s)-P_{LR}(s)$.

The following results are obtained based on the assumption that $\Delta_1>0,\Delta_2<0$ and $\vert\Delta_1\vert<\vert\Delta_2\vert$.

\begin{lemma}\label{LemmaRhombusEstimation}Consider $P_{RL}(s)$ and $P_{LR}(s)$ for different $s$, we have

(1) if $ s<\frac{\Delta_1+\Delta_2}{2}$, then $P_{RL}(s)>P_{LR}(s)$;

(2) if $s>0$, then $P_{RL}(s)<P_{LR}(s)$.

% \begin{enumerate} 
%     \item  if $ s<\frac{\Delta_1+\Delta_2}{2}$, then $P_{RL}(s)>P_{LR}(s)$;
%     \item  if $s>0$, then $P_{RL}(s)<P_{LR}(s)$.
% \end{enumerate}
\end{lemma}

\begin{proof}
Let $g_1(x)= \frac{1-e^x}{x}$. By simple calculation, we can prove that $g_1^\prime(x) \leq 0$, i.e., $g_1(x)$ decreases with $x$. Since $\Delta_2<0<\Delta_1$, we have
\begin{equation*}
    g_1(\eta\Delta_2)=\frac{1-e^{\eta \Delta_2}}{\eta \Delta_2}>\frac{1-e^{\eta \Delta_1}}{\eta \Delta_1}=g_1(\eta\Delta_1),
\end{equation*}
which equals
\begin{equation*}
    f_1 = \Delta_1(1-e^{\eta \Delta_2})-\Delta_2 (1-e^{\eta \Delta_1})<0.
\end{equation*}
Thus, $f(s)$ decreases with $s$. 

Now, we compute
\begin{equation*}
    \begin{aligned}
       f(0)=  -\Delta_2(1-e^{\eta \Delta_1})(1+e^{\eta \Delta_2})+\Delta_1(1-e^{\eta \Delta_2})(1+e^{\eta \Delta_1}).
    \end{aligned}
\end{equation*}

Let $g_2(x)= \frac{e^x-1}{x(1+e^x)}$. Then, by calculation, $g_2(-x)=g_2(x)$ and $g_2(x)$ is monotonically decreasing when $x>0$. Since $\vert\Delta_1\vert<\vert\Delta_2\vert$, we have
\begin{equation*}
    \begin{aligned}
    & g_2(\eta\Delta_1) = \frac{e^{\eta \Delta_1}-1}{\eta \Delta_1(1+e^{\eta \Delta_1})}>\frac{e^{\eta \Delta_2}-1}{\eta \Delta_2(1+e^{\eta \Delta_2})}=g_2(\eta\Delta_2) \\
    \Leftrightarrow\quad & {  -\Delta_2(1-e^{\eta \Delta_1})(1+e^{\eta \Delta_2})+\Delta_1(1-e^{\eta \Delta_2})(1+e^{\eta \Delta_1})   }<0,
    \end{aligned}
\end{equation*}
which means $f(0)<0$. Since $f(s)$ decreases with $s$, we have $f(s)<0$ for all $s>0$. Then, $P_{LR}(s)>P_{RL}(s)$ for all $s>0$.

On the other hand, 
\begin{equation*}
    \begin{aligned}
       f\left(\frac{\Delta_1+\Delta_2}{2}\right)=(e^{\eta\Delta_1/{2}}+e^{\eta\Delta_2/{2}})\left(-\Delta_2e^{\eta\Delta_2/{2}}(1-e^{\eta\Delta_1})+\Delta_1e^{\eta\Delta_1/{2}}(1-e^{\eta\Delta_2})\right).
    \end{aligned}
\end{equation*}

Let $g_3(x)=\frac{xe^{x/2}}{1-e^x}$. Then, by calculation, we have $g_3(-x)=g_3(x)$ and $g_3(x)$ increases with $x$ when $x>0$. Since $\vert\Delta_1\vert<\vert\Delta_2\vert$, we have
\begin{equation*}
    \begin{aligned}
    &g_3(\eta\Delta_2)=\frac{\eta \Delta_2 e^{\eta\Delta_2/{2}}}{1-e^{\eta\Delta_2}}>\frac{\eta \Delta_1 e^{\eta\Delta_1/{2}}}{1-e^{\eta\Delta_1}}=g_3(\eta\Delta_1)\\
    \Leftrightarrow\quad & -\Delta_2e^{\eta\Delta_2/{2}}(1-e^{\eta\Delta_1})+\Delta_1e^{\eta\Delta_1/{2}}(1-e^{\eta\Delta_2})>0
    \end{aligned}
\end{equation*}
which equals $f\left(\frac{\Delta_1+\Delta_2}{2}\right)>0$. Since $f(s)$ decreases with $s$, we have $f(s)>0$ for all $s<\frac{\Delta_1+\Delta_2}{2}$. Hence $P_{LR}(s)<P_{RL}(s)$ for all $s<\frac{\Delta_1+\Delta_2}{2}$.
\end{proof}

When $\Delta_1 = -\Delta_2$, by \ref{LemmaRhombusEstimation}, for any state $s$, we can always compare $P_{LR}(s)$ and $P_{RL}(s)$. Thus, the analysis and results are all valid. 

Denote $s^{0\ast}$ to be the zero point of $f(s)$ defined by equation \eqref{f(s)Def}, i.e., $f(s^{0\ast})=0$. By Lemma \ref{LemmaRhombusEstimation}, we know that 
\begin{equation*}
    \frac{\Delta_1+\Delta_2}{2}<s^{0\ast}<0.
\end{equation*}
In the proof of Lemma \ref{LemmaRhombusEstimation}, we have proven that $f(s)$ decreases with $s$. 
Then, we know
\begin{equation}\label{s0Difference}
    \left\{ \begin{array}{lc}
    P_{LR}(s)<P_{RL}(s),&\quad  \text{if}\ s<s^{0\ast};   \\
    P_{LR}(s)\geq P_{RL}(s),&\quad \text{if}\ s\geq s^{0\ast}. 
    \end{array}
    \right.
\end{equation}

\begin{lemma}\label{Lemmas0Difference} 
Suppose that the optimal action of some states at time $t+1$ is already known. Consider the optimal action of state $s_{i,t}$ at time $t$, we have

(1) if $s_{i,t}<s^{0\ast}$ and $y^\ast(s_{i,t+1})=e_R$, then  $y^\ast(s_{i,t})=e_R$;

(2) if $s_{i,t}\geq s^{0\ast}$ and $y^\ast(s_{i+1,t+1})=e_L$, then $y^\ast(s_{i,t})=e_L$.
% \begin{enumerate}
%     \item if $s_{i,t}<s^{0\ast}$ and $y^\ast(s_{i,t+1})=e_R$, then  $y^\ast(s_{i,t})=e_R$;
%     \item  if $s_{i,t}\geq s^{0\ast}$ and $y^\ast(s_{i+1,t+1})=e_L$, then $y^\ast(s_{i,t})=e_L$.
% \end{enumerate}
\end{lemma}

\begin{proof}
If $s_{i,t}<s^{0\ast}$ and $y^\ast(s_{i,t+1})=e_R$, by \eqref{StateTransFormula},
we have 
\begin{equation*}
    \begin{aligned}
    r(s_{i,t},e_R)+f^\ast(h(s_{i,t},e_R))& = r(s_{i,t},e_R)+f^\ast(s_{i+1,t+1})\\
    &\geq r(s_{i,t},e_R)+r(s_{i+1,t+1},e_L)+f^\ast(s_{i+1,t+2})\\
    & > r(s_{i,t},e_L)+r(s_{i,t+1},e_R)+f^\ast(s_{i+1,t+2})\\
    & = r(s_{i,t},e_L)+f^\ast(s_{i,t+1})\\
    & = r(s_{i,t},e_L)+f^\ast(h(s_{i,t},e_L)),
    \end{aligned}
\end{equation*}
where the first inequality follows from the definition of $f^\ast(s_{i+1,t+1})$, the second inequality follows from formula \eqref{s0Difference} and the second equality holds because $y^\ast(s_{i,t+1})=e_R$. Recall that
$$y^\ast(s_{i,t}) = \arg\max\limits_{y\in \{e_L,e_R\}}\{r(s_{i,t},y)+f^\ast(h(s_{i,t},y))\}$$
hence $y^\ast(s_{i,t})=e_R$.

Similarly, if $s_{i,t}\geq s^{0\ast}$ and $y^\ast(s_{i+1,t+1})=e_L$, 
we have $y^\ast(s_{i,t})=e_L$.

% \begin{equation*}
%     \begin{aligned}
%     r(s_{i,t},e_L)+f^\ast(h(s_{i,t},e_L))& =
%     r(s_{i,t},e_L)+f^\ast(s_{i,t+1})\\
%     &\geq r(s_{i,t},e_L)+r(s_{i,t+1},e_R)+f^\ast(s_{i+1,t+2})\\
%     & > r(s_{i,t},e_R)+r(s_{i+1,t+1},e_L)+f^\ast(s_{i+1,t+2})\\
%     & = r(s_{i,t},e_R)+f^\ast(s_{i,t+1})\\
%     & = r(s_{i,t},e_R)+f^\ast(h(s_{i,t},e_R))
%     \end{aligned}
% \end{equation*}
% where the first inequality follows by the definition of $f^\ast(s_{i+1,t+1})$ and the second inequality follows by formula \eqref{s0Difference}. Therefore, 
\end{proof}

Since $-\vert\Delta_2\vert<\frac{\Delta_1+\Delta_2}{2}<s^{0\ast}<0,$ we can easily obtain the following corollary.
\begin{corollary}\label{LemmaRtoRLtoL}
Consider the optimal action of state $s_{i,t}$ at time $t$, we have

(1) if $s_{i,t}<-\vert\Delta_2\vert$ and $y^\ast(s_{i,t+1})=e_R$, then  $y^\ast(s_{i,t})=e_R$;

(2) if $s_{i,t}>0$ and $y^\ast(s_{i+1,t+1})=e_L$, then $y^\ast(s_{i,t})=e_L$.

\end{corollary}

Corollary \ref{LemmaRtoRLtoL} gives the recurrence relation between the optimal actions of states, which is crucially important in the following study of the optimal play of player Y.

\subsection{Property of the Optimal Path}

Now, we will study the property of the optimal play of player Y and the optimal path of the system. We present and prove the results for the case when $s^\ast>0$. The result when $s^\ast\leq 0$ can be obtained in the same way.

To obtain the optimal play, we need to get the optimal action of every state by solving the Bellman Optimality Equation
\begin{equation*}
   y^\ast(s_t) = \arg\max\limits_{y\in \{e_L,e_R\}}\{r(s_{t},y)+f^\ast(h(s_t,y))\}.
\end{equation*}

Clearly, it is costly and infeasible to directly solve this equation for every state. However, by using Corollary \ref{LemmaRtoRLtoL}, the optimal action of states could be recurrently obtained and costly computation could be avoided. Moreover, the periodic performance of the zero path could further reduce the computation to obtain the optimal play, which will be shown below.

First, consider the system states at time $T$. Denote $j^\ast$ to be the index such that
\begin{equation}\label{jast}
    s_{j^\ast,T}\leq s^\ast,\quad s_{j^\ast+1,T}>s^\ast.
\end{equation}
Since $T$ is sufficiently large, we have $s^\ast>\vert\Delta_2\vert$ and then $s_{j^\ast,T}> \vert\Delta_2\vert$.

Since $f^\ast(s_{T+1})\equiv 0$, by \eqref{yast} and by \eqref{MBR}, we have 
\begin{equation}\label{yastT}
    y^\ast(s_{i,T})=\begin{cases}
    e_R, & i\leq j^\ast;\\
    e_L, & i\geq j^\ast+1.
    \end{cases}
\end{equation}

Recall that sequence $s_1^0,s_2^0,\cdots,s_{T+1}^0$ is called the zero sequence and sequence $s_1^0,y_1^0,s_2^0,y_2^0,\cdots,s_T^0,y_T^0,s_{T+1}^0$ is called the zero path.

Define $i^0(t)$ to be the index such that $s_{i^0(t),t}=s_t^0$. By \eqref{StateTransFormula} and \eqref{MBR}, we have
\begin{equation}\label{i0Relation}
    i^0(t) = \begin{cases}
    i^0(t-1), &\quad \text{if}\ s_{t-1}^0\geq 0;\\
    i^0(t-1)+1, & \quad \text{if}\ s_{t-1}^0<0.
    \end{cases}
\end{equation}

From state $s_{j^\ast,T}$, we consider a state sequence $s_{j^\ast,T},s_{j^\ast-1,T-1}, s_{j^\ast-2,T-2},\cdots$ backwards. Note that the sequence satisfies $s_{j^\ast-k+1,T-k+1}=h(s_{j^\ast-k,T-k},e_R), k\geq 1.$ Thus, these states can be connected by a straight line on STTG, called \textit{the exploring path}.

\begin{lemma}\label{Lemma1}
The zero path and the exploring path have an intersection point at some time $t^{\times}$ satisfying $T-t^{\times} = o(T)$ i.e., there exists time $t^{\times}$ such that
$$s_{t^{\times}}^0 = s_{j^\ast,T}-(T-t^{\times})\vert\Delta_2\vert=s_{j^\ast-T+t^{\times},t^{\times}}.$$
Further, for time $t$, $t^{\times}\leq t\leq T$, we have
\begin{equation}\label{NeedToExplain}
    y^\ast(s_{i,t})  = \begin{cases}
    e_L,& i\geq j^\ast-T+t+1;\\
    e_R, & i\leq i^0(t)-2.
    \end{cases}
\end{equation}
\end{lemma}

\begin{figure}[htbp]
    \centering
    \includegraphics[width=0.9\textwidth]{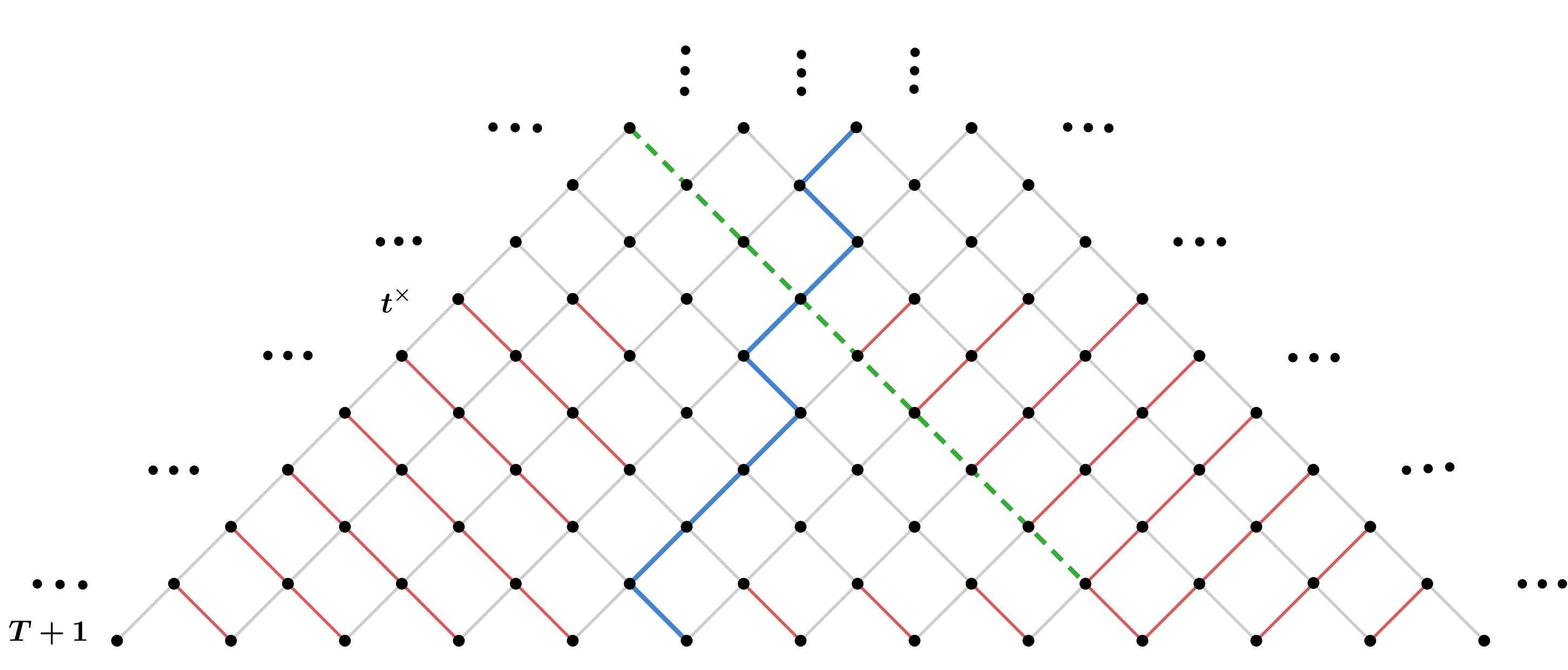}
    \caption{Illustration of Lemma \ref{Lemma1}. }
    \label{fig:intersection}
\end{figure}

The result of Lemma \ref{Lemma1} is demonstrated by Figure \ref{fig:intersection}. The blue polyline represents the zero path and the green broken line represents the exploring path. As shown in the graph, the zero path and the exploring path have an interaction point at time $t^{\times}$. For states in the index range given by \eqref{NeedToExplain}, their optimal actions can be obtained, which are represented by red line segments. 

Equation \eqref{NeedToExplain} actually gives a range of states whose optimal action can be obtained. 
As shown in Figure \ref{fig:intersection}, on the one hand, the optimal action of states at the right side of the exploring path can be computed; on the other hand, for state $s_{i,t}$ at the left side of the myopic path, if the absolute value of the difference between its index $i$ and $i^0(t)$ is no less than $2$, its optimal action $y^\ast(s_{i,t})$ can also be computed. 

Now, we prove Lemma \ref{Lemma1}.
\begin{proof}
First, consider the state sequence $s_{j^\ast,T},s_{j^\ast-1,T-1},\cdots$ on the exploring path, we have
\begin{equation}\label{TEMP1}
    s_{j^\ast-T+t,t} = s_{j^\ast,T}-(T-t)\vert\Delta_2\vert.
\end{equation}
Since $s_{j^\ast,T}>0$, there exists a state $s_{j^\ast-T+t^{\times}-1,t^{\times}-1}$ satisfying $s_{j^\ast-T+t^{\times}-1,t^{\times}-1}<0$, $s_{j^\ast-T+t^{\times},t^{\times}-1}\geq 0$ and $s_{j^\ast-T+t^{\times},t^{\times}}\geq 0$. By Lemma \ref{LemmaProperties}, we have $s_{j^\ast-T+t^{\times},t^{\times}}=s_{t+1}^0$. 
This means that the exploring path has an intersection point with the zero path at time $t^{\times}$.  

By \eqref{jast}, $s_{j^\ast,T}\leq s^\ast$. By \eqref{TEMP1}, since $s_{j^\ast-T+t^{\times},t^{\times}}\geq 0$, we have 
\begin{equation}\label{Ranget0}
    T-t^{\times}\leq \frac{s_{j^\ast,T}}{\vert\Delta_2\vert}\leq \frac{s^\ast}{\vert\Delta_2\vert}\leq \lceil \frac{s^\ast}{\vert\Delta_2\vert}\rceil.
\end{equation}
Recall that $s^\ast = o(T)$, then $T-t^{\times} = o(T)$.

Then, we prove that for state $s_{i,t}$ where $t^{\times}\leq t\leq T$ and $i\geq j^\ast-T+t+1$, we have $y^\ast(s_{i,t})=e_L$. 
We use mathematical induction to prove this. 

When $t=T$, by \eqref{yastT}, for state $s_{i,T}$ where $i\geq j^\ast+1$, we have $y^\ast(s_{i,T})=e_L$. Assume that for state $s_{i,t+1}$ where $t^{\times}+1\leq t+1\leq T$ and $i\geq j^\ast-T+t+2$, we have $y^\ast(s_{i,t+1})=e_L$. Then, for state $s_{i,t}$ where $i\geq j^\ast-T+t+1$,  $$s_{i,t}\geq s_{j^\ast-T+t+1,t}\geq s_{j^\ast-T+t^{\times}+1,t^{\times}}>s_{j^\ast-T+t^{\times},t^{\times}-1}\geq 0.$$ By induction assumption that $y^\ast(s_{i+1,t+1})=e_L$ and by Corollary \ref{LemmaRtoRLtoL}, we have $y^\ast(s_{i,t})=e_L.$

Lastly, we prove that for state $s_{i,t}$ where $t^{\times}\leq t\leq T$ and $i\leq i^0(t)-2$, we have $y^\ast(s_{i,t})=e_R$. 
We also use mathematical induction. 

When $t=T$, by \eqref{yastT}, for state $s_{i,T}$ where $i\leq i^0(T)-2$, we have $y^\ast(s_{i,T})=e_R$. Assume that for state $s_{i,t+1}$ where $t^{\times}+1\leq t+1\leq T$ and $i\leq i^0(t+1)-2$, we have $y^\ast(s_{i,t+1})=e_R$. Then, for state $s_{i,t}$ where $i\leq i^0(t)-2$, $$s_{i,t}\leq s_{i^0(t)-2,t}=s_t^0-2(\vert\Delta_2\vert+\vert\Delta_1\vert) < \vert\Delta_2\vert-2(\vert\Delta_2\vert+\vert\Delta_1\vert)< -\vert\Delta_2\vert.$$ By \eqref{i0Relation}, $i^0(t)\leq i^0(t+1)$ and then from the induction assumption, $y^\ast(s_{i,t+1})=e_R$ for $i\leq i^0(t)-2\leq i^0(t+1)-2$. Therefore, by Corollary \ref{LemmaRtoRLtoL}, we have $y^\ast(s_{i,t})=e_L.$
\end{proof}

From time $t^{\times}$ backwards, following the zero path, from Lemma \ref{LemmaProperties}, we know that there exists time $t^d\leq t^{\times}$ such that
\begin{equation}\label{tddef}
 s_{t^d}^0<0,s_{t^d+1}^0\geq 0.   
\end{equation}
Recall that the least positive integer of the zero sequence is $T^\ast = \frac{m(\vert\Delta_1\vert,\vert\Delta_2\vert)}{\vert\Delta_2\vert}+\frac{m(\vert\Delta_1\vert,\vert\Delta_1\vert)}{\vert\Delta_2\vert}$. By \eqref{Ranget0}, we have
\begin{equation}\label{Rangetd}
    T-t^d\leq T-t^{\times}+T^\ast\leq \lceil \frac{s^\ast}{\vert\Delta_2\vert}\rceil + T^\ast.
\end{equation}
Since $s^\ast = o(T)$, we know that $T-t^d = o(T)$.

Then, consider the optimal actions of states at time $t^d-1$.

\begin{lemma}\label{Lemma2} For all of states at time $t^d-1$, their optimal actions  
can be obtained. To be more specific, we have
\begin{equation*}
    y^\ast(s_{i,t^d-1})=\begin{cases}
    e_R, & i\leq i^0(t^d-1)-1;\\
    e_L, & i\geq i^0(t^d-1).
    \end{cases}
\end{equation*}
\end{lemma}

\begin{figure}[htbp]
    \centering
    \includegraphics[width=0.90\textwidth]{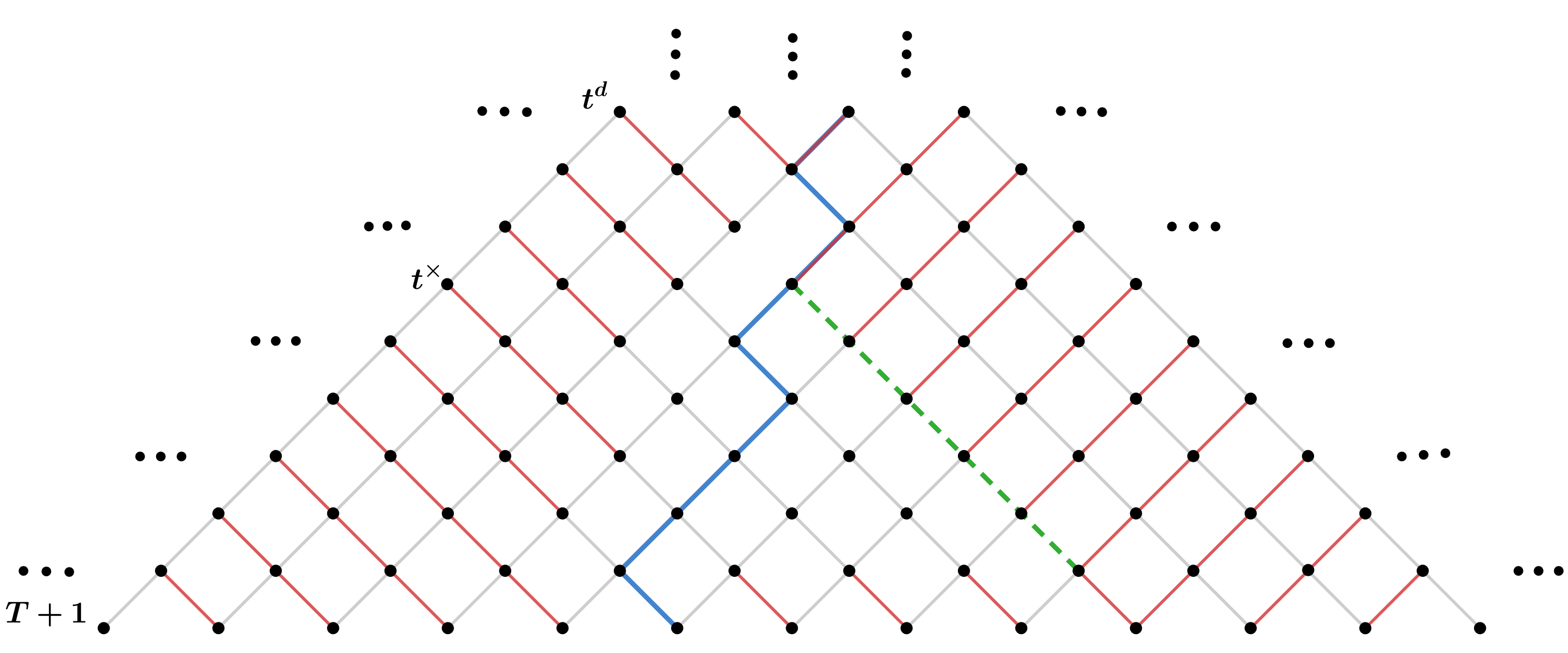}
    \caption{Illustration of Lemma \ref{Lemma2}. In the graph, state $s_{t^d}^0$ on the zero path is less than 0 and $s_{t^d+1}^0\geq 0$. For states at time $t^d-1$, their optimal actions can all be obtained.}
    \label{fig:td}
\end{figure}

Figure \ref{fig:td} illustrates the result of Lemma \ref{Lemma2} and can help understand the proof of Lemma \ref{Lemma2}.

\begin{proof}
By Lemma \ref{Lemma1}, for states at time $t^{\times}$, we have
\begin{equation}\label{Lemma5.4Eq1}
    y^\ast(s_{i,t^{\times}})=\begin{cases}
    e_L, & i\geq i^0(t^{\times})+1;\\
    e_R, & i\leq i^0(t^{\times})-2.
    \end{cases}
\end{equation}

We prove Lemma \ref{Lemma2} for different cases.

\textbf{{Case 1}}: $s_{t^{\times}-1}^0<0$. In this case, by \eqref{i0Relation}, $i^0(t^{\times}-1)=i^0(t^{\times})-1$. For state $s_{i,t^{\times}-1}$ where $i\geq i^0(t^{\times}-1)+1$, by Lemma \ref{LemmaProperties}, we have $$s_{i,t^{\times}-1}\geq s_{i^0(t^{\times}-1)+1,t^{\times}-1} =  s_{t^{\times}-1}^0+\vert\Delta_1\vert+\vert\Delta_2\vert\geq -\vert\Delta_1\vert+\vert\Delta_1\vert+\vert\Delta_2\vert>0.$$ 
For $i\geq i^0(t^{\times}-1)+1$, since $i^0(t^{\times}-1)=i^0(t^{\times})-1$, we have $i+1\geq i^0(t)+1$ and then by \eqref{Lemma5.4Eq1}, 
$y^\ast(s_{i+1,t^{\times}})=e_L$. Thus, by Corollary \ref{LemmaRtoRLtoL}, we have $y^\ast(s_{i,t^{\times}-1})=e_L$ for $i\geq i^0(t^{\times}-1)+1$. Similarly, we can prove that $y^\ast(s_{i,t^{\times}-1})=e_R$ for $i\leq i^0(t^{\times}-1)-1$. For the states at time $t^{\times}-1$, we have
\begin{equation}\label{yastt0-1}
    y^\ast(s_{i,t^{\times}-1})=\begin{cases}
    e_L, & i\geq i^0(t^{\times}-1)+1;\\
    e_R, & i\leq i^0(t^{\times}-1)-1.
    \end{cases}
\end{equation}

Then, consider the states at time $t^{\times}-2$. By Lemma \ref{LemmaProperties}, since $s_{t^{\times}-1}^0<0$, we have $s_{t^{\times}-2}^0\geq 0$ and then by \eqref{i0Relation}, $i^0(t^{\times}-2)=i^0(t^{\times}-1)$. For state $s_{i,t^{\times}-2}$ where $i\geq i^0(t^{\times}-2)$,
$$s_{i,t^{\times}-2}\geq s_{t^{\times}-2}^0\geq 0.$$  
For $i\geq i^0(t^{\times}-2)$, since $i^0(t^{\times}-2)=i^0(t^{\times}-1)$, we have $i+1\geq i^0(t^{\times}-1)+1$ and then by \eqref{yastt0-1}, 
$y^\ast(s_{i+1,t^{\times}-1})=e_L$. Thus, by Corollary \ref{LemmaRtoRLtoL}, we have $$y^\ast(s_{i,t^{\times}-2})=e_L$$ for $i\geq i^0(t^{\times}-2)$.

On the other hand, for state $s_{i,t^{\times}-2}$ where $i\leq i^0(t^{\times}-2)-1$, since $i^0(t^{\times}-2)=i^0(t^{\times}-1)$,
we have 
\begin{equation*}
    s_{i,t^{\times}-2}\leq s_{i^0(t^{\times}-2)-1,t^{\times}-2}=s_{i^0(t^{\times}-2),t^{\times}-1}-\vert\Delta_2\vert=s_{t^{\times}-1}^0-\vert\Delta_2\vert<-\vert\Delta_2\vert.
\end{equation*}
By \eqref{yastt0-1}, for $i\leq i^0(t^{\times}-2)-1=i^0(t^{\times}-1)-1$, we have $y^\ast(s_{i,t^{\times}-1})=e_R$. Thus, by Corollary \ref{LemmaRtoRLtoL}, we have $$y^\ast(s_{i,t^{\times}-2})=e_R$$ for $i\leq i^0(t^{\times}-2)-1.$

Up to now, we have proven that 
\begin{equation*}
    y^\ast(s_{i,t^{\times}-2}) = \begin{cases}
    e_L, & i\geq i^0(t^{\times}-2);\\
    e_R, & i\leq i^0(t^{\times}-2)-1.
    \end{cases}
\end{equation*}
Then, we find $t^d=t^{\times}-1$. 

\textbf{{Case 2:}} $s_{t^{\times}-1}^0\geq 0$. By \eqref{i0Relation}, $i^0(t^{\times}-1)=i^0(t^{\times})$. Consider the optimal action of states at time $t^{\times}-1$. For state $s_{i,t^{\times}-1}$ where $i\geq i^0(t^{\times}-1)$, on one hand, $$s_{i,t^{\times}-1}\geq s_{i^0(t^{\times}-1),t^{\times}-1}=s_{t^{\times}-1}^0\geq 0;$$ on the other hand, for $i\geq i^0(t^{\times}-1)$, since $i^0(t^{\times}-1)=i^0(t^{\times})$, we have $i+1\geq i^0(t^{\times})+1$ and then by \eqref{Lemma5.4Eq1}, $y^\ast(s_{i+1,t^{\times}})=e_L.$ Thus, by Corollary \ref{LemmaRtoRLtoL}, we have 
\begin{equation*}
    y^\ast(s_{i,t^{\times}-1})=e_L
\end{equation*}
for $i\geq i^0(t^{\times}-1).$
Similarly, we have $$y^\ast(s_{i,t^{\times}-1})=e_R$$ 
for $i\leq i^0(t^{\times}-1)-1.$ 

In conclusion, when $s_{t^{\times}-1}\geq 0$, for states at time $t^{\times}-1$, we have proven that 
\begin{equation*}
    y^\ast(s_{i,t^{\times}-1}) = 
    \begin{cases}
    e_L, & i\geq i^0(t^{\times}-1);\\
    e_R, & i\leq i^0(t^{\times}-1)-1.
    \end{cases}
\end{equation*}

Then, consider states at time $t^{\times}-2$. 
If $s_{t^{\times}-2}^0<0$, this case is the same with \textbf{case 1} and thus Lemma \ref{Lemma2} holds. 
If $s_{t^{\times}-2}\geq 0$, we can prove
\begin{equation*}
    y^\ast(s_{i,t^{\times}-2}) = 
    \begin{cases}
    e_L, & i\geq i^0(t^{\times}-2);\\
    e_R, & i\leq i^0(t^{\times}-2)-1;
    \end{cases}
\end{equation*}
in the same way as the case that $s_{t^{\times}-1}^0\geq 0$. In this case, $s_{t^{\times}-2}^0 = s_{t^{\times}}^0+2\vert \Delta_1 \vert$. 

By Lemma \ref{LemmaProperties}, we know that $s_t^0\leq \vert \Delta_2 \vert$ for all time $t$. Therefore, we will encounter 
\textbf{Case 1} at some time $t\leq t^{\times}-2$ and then Lemma \ref{Lemma2} holds.
\end{proof}

By now, we have computed the optimal actions for states at time $t\geq t^d-1$ as many as possible. Next, we consider the states at time $t\leq t^d-1$. 
% Actually, if we can compute the optimal action of every state at time $t\leq t^d-1$, the optimal play from time $1$ to time $t^d-1$ is then obtained.

\begin{lemma}\label{Lemma3} For state $s_{i,t}$ where $t\leq t^d-1$,
if $s_t^0\geq 0$, we have
\begin{equation*}
   y^\ast(s_{i,t}) = \begin{cases}
   e_L, & i\geq i^0(t);\\
   e_R, & i\leq i^0(t)-2;
   \end{cases}
\end{equation*}
while if $s_t^0<0$, we have 
\begin{equation*}
    y^\ast(s_{i,t}) = \begin{cases}
    e_L, & i\geq i^0(t)+1;\\
    e_R, & i\leq i^0(t)-1.
    \end{cases}
\end{equation*}
\end{lemma}

\begin{figure}[htbp]
    \centering
    \includegraphics[width=0.90\textwidth]{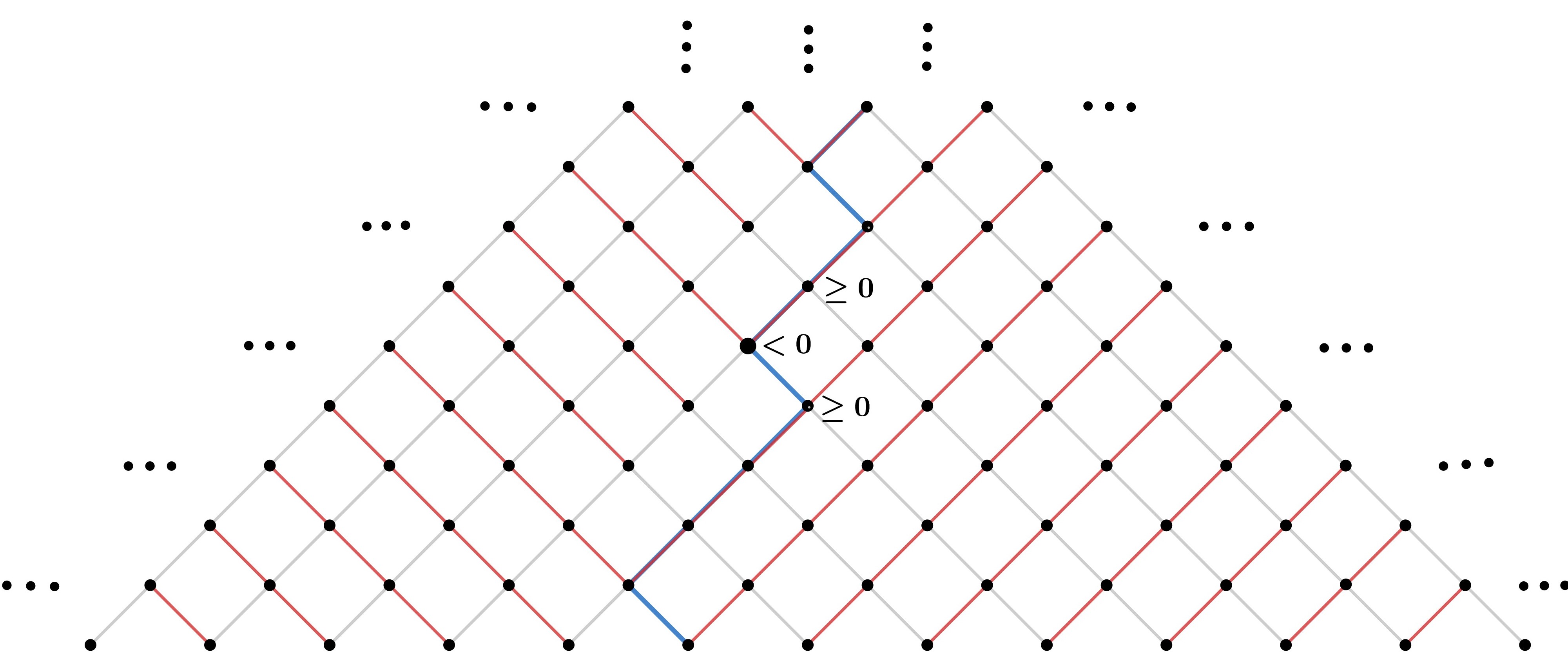}
    \caption{Illustration of Lemma \ref{Lemma3}. From the graph, we can see that if $s_t^0<0$, then only $y^\ast(s_{i^0(t),t})$ cannot be determined; if $s_t^0\geq 0,$ then only $y^\ast(s_{i^0(t)-1,t})$ cannot be determined.}
    \label{fig:yast}
\end{figure}

Figure \ref{fig:yast} illustrates the result of Lemma \ref{Lemma3} and can help understand the proof of Lemma \ref{Lemma3}.

\begin{proof}
We prove Lemma \ref{Lemma3} by mathematical induction. 
First, by Lemma \ref{Lemma2}, we know that Lemma \ref{Lemma3} holds for time $t = t^d-1$. 

Assume that Lemma \ref{Lemma3} holds for time $t\leq t^d-1$ and then consider the optimal action of the states at time $t-1$. 

We discuss this problem for different cases.

\textbf{Case 1:} $s_t^0\geq 0$. For state $s_{i,t-1}$ where $i\geq i^0(t)$, we have $$s_{i,t-1}\geq s_{i^0(t),t-1}=s_{i^0(t),t}+\vert\Delta_1\vert = s_t^0+\vert\Delta_1\vert>0.$$ 
Since we have assumed that Lemma \ref{Lemma3} holds for time $t$ and $s_t^0\geq 0$, then $y^\ast(s_{i,t})=e_L$ for 
$i\geq i^0(t)$. 
Similarly, we can prove that $y^\ast(s_{i,t-1})=e_R$ for $i\leq i^0(t)-2$. Therefore, 
\begin{equation}\label{Lemma5.5Eq1}
    y^\ast(s_{i,t-1}) = 
    \begin{cases}
    e_L, & i\geq i^0(t);\\
    e_R, & i\leq i^0(t)-2.
    \end{cases}
\end{equation}

Further, consider the sign of $s_{t-1}^0$:

\textbf{Case 1.1: } $s_{t-1}^0<0$. By \eqref{i0Relation}, $i^0(t-1)=i^0(t)-1$. Thus, formula \eqref{Lemma5.5Eq1} is equivalent to 
    \begin{equation*}
    y^\ast(s_{i,t-1}) = \begin{cases}
    e_L, & i\geq i^0(t-1)+1\\
    e_R, & i\leq i^0(t-1)-1.
    \end{cases}
    \end{equation*}
    
    \textbf{Case 1.2:} $s_{t-1}^0\geq 0$. By \eqref{i0Relation}, $i^0(t-1)=i^0(t)$. Thus, formula \eqref{Lemma5.5Eq1} is equivalent to
    \begin{equation*}
   y^\ast(s_{i,t-1}) = 
   \begin{cases}
   e_L, & i\geq i^0(t-1);\\
   e_R, & i\leq i^0(t-1)-2.\\
   \end{cases}
\end{equation*}

\textbf{Case 2:} $s_t^0<0$. By Lemma \ref{LemmaProperties},  $s_{t-1}^0\geq 0$ and then by \eqref{i0Relation}, $i^0(t)=i^0(t-1)$.
For state $s_{i,t-1}$ where $i\leq i^0(t-1)-1$, we have 
\begin{equation*}
    s_{i,t-1}\leq s_{i^0(t-1)-1,t-1}=s_{i^0(t),t}-\vert\Delta_2\vert=s_t^0-\vert\Delta_2\vert<-\vert\Delta_2\vert.
\end{equation*}
For $i\leq i^0(t-1)-1$, since $i^0(t)=i^0(t-1)$, we have $i\leq i^0(t)-1$ and then from the induction assumption, since 
$s_t^0<0$, we have $y^\ast(s_{i,t})=e_R$ for $i\leq i^0(t)-1$. Thus, by Corollary \ref{LemmaRtoRLtoL}, 
\begin{equation*}
    y^\ast(s_{i,t-1})=e_R
\end{equation*}
for $i\leq i^0(t-1)-1$. 

On the other hand, for state $s_{i,t-1}$ where $i\geq i^0(t-1)$, we have 
\begin{equation*}
  s_{i,t-1}\geq s_{i^0(t-1),t-1}\geq 0. 
\end{equation*}
For $i\geq i^0(t-1)$, since $i^0(t)=i^0(t-1)$, we have $i\geq i^0(t)$ and then from the induction assumption, since 
$s_t^0<0$, we have $y^\ast(s_{i,t})=e_L$ for $i\geq i^0(t)$. Thus, by Corollary \ref{LemmaRtoRLtoL}, 
\begin{equation*}
    y^\ast(s_{i,t-1})=e_L
\end{equation*}
for $i\geq i^0(t-1)$. 

In conclusion, if $s_t^0<0$, we have $s_{t-1}^0\geq 0$ and 
\begin{equation}\label{Lemma5.5Eq2}
    y^\ast(s_{i,t-1}) = 
   \begin{cases}
   e_L, & i\geq i^0(t-1);\\
   e_R, & i\leq i^0(t-1)-1.\\
   \end{cases}
\end{equation}

So far, we have discussed all possible cases of $s_{t}^0$ and $s_{t-1}^0$. For every case, we have proven Lemma \ref{Lemma3} holds for time $t-1$.  Therefore, by the method of mathematical induction, Lemma \ref{Lemma3} is correct for all time $t\leq t^d-1$.
\end{proof}

Lemma \ref{Lemma3} means that for time $t\leq t^d-1$, there is only one state at time $t$ whose optimal action has not yet been computed. 
Even so, the optimal play of player Y cannot be determined. 

From \textbf{Case 2} in the proof of Lemma \ref{Lemma3}, we actually have proven the corollary below.

\begin{corollary}\label{Corollary1}
If time $t$ satisfies $s_{t}^0\geq 0$ and $s_{t+1}^0<0$, then all of the optimal actions of states at time $t$ can be obtained. To be more specific, we have
\begin{equation}\label{Lemma5.5Eq3}
    y^\ast(s_{i,t}) = \begin{cases}
    e_L, & i\geq i^0(t);\\
    e_R, & i\leq i^0(t)-1.
    \end{cases}
\end{equation}
\end{corollary}

Based on Corollary \ref{Corollary1} and Lemma \ref{Lemma3}, we will show that state $s_t^\ast$, which is on the optimal path, could be obtained for some time $t$.

\begin{lemma}\label{Lemma4} For state $s_t^0$ at time $t\leq t^d$, 

(i) if $s_t^0<0$, then $s_t^\ast=s_t^0$;

(ii) if $s_t^0\geq 0$, then $s_t^\ast=s_t^0$ or $s_t^\ast=s_{i^0(t)-1,t}$ where $s_{i^0(t)-1,t}$ is the state next to $s_t^0$ on its left side on STTG.
\end{lemma}

\begin{proof}
We prove Lemma \ref{Lemma4} using mathematical induction. First, for $t=1$, $s_1^\ast=s_1^0=0$ and by \eqref{MBR}, $\operatorname{MBR}(s_1^0)=e_L$. Then, for $t=2$, $s_2^0=-\vert\Delta_1\vert<0$. By Lemma \ref{Lemma3}, $y^\ast(s_1^\ast)=y^\ast(s_1^0)=e_L$ and then we have $s_2^\ast = s_2^0$. 

Assume that Lemma \ref{Lemma4} holds for time $t$, $1\leq t\leq t^d)$ and then consider $s_{t+1}^\ast$ for time $t+1$ by classified discussion.

\textbf{Case 1:} $s_t^0<0$. By Lemma \ref{LemmaProperties}, $s_{t+1}^0>0$ and then, by \eqref{i0Relation}, $i^0(t+1)=i^0(t)+1$. From the induction assumption, we have $s_t^\ast=s_t^0$ then either $s_{t+1}^\ast = s_{i^0(t),t+1}=s_{i^0(t+1)-1,t+1}$ or $s_{t+1}^\ast = s_{i^0(t)+1,t+1}=s_{i^0(t+1),t+1}$. 

\textbf{Case 2:} $s_t^0\geq 0$ and $s_{t+1}^0\geq 0$. Then, by \eqref{i0Relation}, $i^0(t+1)=i^0(t)$. From the induction assumption, if $s_t^\ast=s_{i^0(t)-1,t}$ then either $s_{t+1}^\ast=s_{i^0(t)-1,t+1}=s_{i^0(t+1)-1,t+1}$ or $s_{t+1}^\ast=s_{i^0(t),t+1}=s_{i^0(t+1),t+1} = s_{t+1}^0$.  If $s_t^\ast=s_{i^0(t),t}$, by Lemma \ref{Lemma3}, $y^\ast(s_{i^0(t),t})=e_L$ and then $s_{t+1}^\ast=s_{i^0(t),t+1}=s_{i^0(t+1),t+1}=s_{t+1}^0$. 
    
    Therefore, if $s_{t+1}^0\geq 0$, then either $s_{t+1}^\ast=s_{i^0(t+1)-1,t+1}$ or $s_{t+1}^\ast=s_{i^0(t+1),t+1}$.
    
 \textbf{Case 3:} $s_t^0\geq 0$ and $s_{t+1}^0< 0$. Then, by \eqref{i0Relation}, $i^0(t+1)=i^0(t)$. By \eqref{Lemma5.5Eq3}, $y^\ast(s_{i^0(t)-1,t})=e_R$ and $y^\ast(s_{i^0(t),t})=e_L$. Recall the induction assumption, either $s_t^\ast=s_{i^0(t)-1,t}$ or $s_t^\ast=s_{i^0(t),t}$. 
    If $s_t^\ast=s_{i^0(t)-1,t}$, then $s_{t+1}^\ast=s_{i^0(t),t+1}=s_{i^0(t+1),t+1}=s_{t+1}^0$; if $s_t^\ast=s_{i^0(t),t}$, then $s_{t+1}^\ast=s_{i^0(t),t+1}=s_{i^0(t+1),t+1}=s_{t+1}^0$. 
    
    Therefore, if $s_{t+1}^0< 0$, then $s_{t+1}^\ast=s_{t+1}^0$.

In conclusion, if $s_{t+1}^0<0$, then $s_{t+1}^\ast=s_{t+1}^0$ 
 while if $s_{t+1}^0\geq 0$, then either $s_{t+1}^\ast=s_{i^0(t+1)-1,t+1}$ {or} $ s_{t+1}^\ast=s_{i^0(t+1),t+1}.$ Thus, we prove Lemma \ref{Lemma4}.
\end{proof}

Lemma \ref{Lemma4} means that if $t\leq t^d$, then all the states $s_t^0$ satisfying $s_t^0<0$ on the zero path are also on the optimal path, i.e., $s_t^\ast = s_t^0$. 

By Lemma \ref{LemmaProperties}, the sequence $s_1^0, s_2^0, \cdots, s_{t^d}^0$ is periodic with period being $T^\ast = \frac{m(\vert\Delta_1\vert,\vert\Delta_2\vert)}{\vert\Delta_2\vert}+\frac{m(\vert\Delta_1\vert,\vert\Delta_2\vert)}{\vert\Delta_1\vert}$. Since $s_1^0 = 0$, by \eqref{MBR}, $\operatorname{MBR}(s_1^0)=e_L$ and then $s_2^0 = s_1^0 - \vert\Delta_1\vert<0$. Hence, $s_2^0$ is on the optimal path, i.e., $s_2^\ast = s_2^0$. Furthermore, $s_{2+T^\ast}^0 = s_{2+2T^\ast}^0 = \cdots = s_2^0<0$ are also on the optimal path, i.e., $s_{2+T^\ast}^\ast = s_{2+T^\ast}^0, s_{2+2T^\ast}^\ast = s_{2+2T^\ast}^0,\cdots$.

Consider the sequence $s_2^\ast, s_{3}^\ast,\cdots, s_{2+T^\ast}^\ast$ and $s_{2+T^\ast}^\ast, s_{3+T^\ast}^\ast,\cdots, s_{2+2T^\ast}^\ast$. 

Recall that finding the optimal path is equivalent to solving the Longest Path Problem on STTG, 
hence the segment of the optimal path from any time $\hat{t}$ to $\tilde{t}>\hat{t}$ is the same with the solution of the Longest Path Problem in \textit{local} STTG (part of the complete STTG from state $s_{\hat{t}}^\ast$ to state $s_{\tilde{t}}^\ast$). Therefore, the sequence $s_2^\ast,s_{3}^\ast,\cdots,s_{2+T^\ast}^\ast$ can be obtained by solving
\begin{equation*}
    \max\limits_{y_2,y_3,\cdots,y_{1+T^\ast}} f^\ast(s_2^\ast,s_{2+T^\ast}^\ast),
\end{equation*}
and the sequence $s_{2+T^\ast}^\ast, s_{3+T^\ast}^\ast,\cdots, s_{2+2T^\ast}^\ast$ can be obtained by solving
\begin{equation*}
    \max\limits_{y_{2+T^\ast},y_{3+T^\ast},\cdots,y_{1+2T^\ast}} f^\ast(s_{2+T^\ast}^\ast,s_{2+2T^\ast}^\ast).
\end{equation*}

Since these two Longest Path Problems have the same starting node, middle nodes and ending node, the local STTGs are actually the same by the generation of STTG, the sequence $s_2^\ast, s_{3}^\ast,\cdots, s_{2+T^\ast}^\ast$ is the same with the sequence $s_{2+T^\ast}^\ast, s_{3+T^\ast}^\ast,\cdots, s_{2+2T^\ast}^\ast$, i.e., for time $t$, $2\leq t\leq 2+T^\ast$
\begin{equation*}
    s_t^\ast = s_{t+T^\ast}^\ast.
\end{equation*}

% For other time $t$ when $s_t^0\geq 0$, by Lemma \ref{Lemma4}, we know that $s_t^\ast\in \{s_{i^0(t)-1,t},s_t^0\}$ where $s_{i^0(t)-1,t}$ is the state next to $s_t^0$ on the left in STTG.

It is easy to see that the above induction holds for the sequence $s_{2+kT^\ast}^\ast,$$ s_{3+kT^\ast}^\ast,$ $\cdots, s_{2+(k+1)T^\ast}^\ast$. Therefore, we can obtain Theorem \ref{ThmBest} finally.

For the state sequence $s_1^\ast, s_2^\ast, \cdots, s_{t^d}^\ast$, we have the following proposition, in which $s^{0\ast}$ is the zero point of function $f(s)$ defined by \eqref{f(s)Def}.
\begin{proposition}\label{lemma_eL_state_greater_than_s0_er_state_less_than_el}
    {For time $1\leq t\leq t^d$, if $y^\ast(s_t^\ast) = e_L$, then $s_t^\ast\geq s^{0\ast}$; if $y^\ast(s_t^\ast) = e_R$, then $s_t^\ast< s^{0\ast}$. }
\end{proposition}

\begin{proof}

    By Lemma \ref{Lemma4}, for time $t\leq t^d$, if $s_t^0<0$, then $s_t^\ast = s_t^0$. Hence, at time $t\leq t^d$ when $s_t^0<0$, the state $s_t^\ast$ is determined and equals $s_t^0$. By Lemma \ref{LemmaProperties}, we can find an interval $[t, t+n), n\geq 2$ such that $s_t^0<0, s_{t+n}^0<0$ and $s_{t+k}^0\geq 0$ for $1\leq k <n$. 

    By Lemma \ref{LemmaProperties}, we obtain $s_{t+n-1}^0\geq 0$. By Corollary \ref{Corollary1}, all of the optimal action at time $t+n-1$ can be obtained, and we have $y^\ast(s_{i^0(t+n-1)-1, t+n-1}) = e_R$ and $y^\ast(s_{i^0(t+n-1), t+n-1}) = e_L$. If $s_{t+n-1}^\ast = s_{i^0(t+n-1)-1, t+n-1}$, then $y^\ast(s_{t+n-1}^\ast) = e_R$ and $s_{t+n-1}^\ast = s_{t+n}^0 -\Delta_2 < -\Delta_2 < s^{0\ast}$; if $s_{t+n-1}^\ast = s_{i^0(t+n-1), t+n-1} = s_{t+n-1}^0$, then $y^\ast(s_{t+n-1}^\ast) = e_L$ and $s_{t+n-1}^\ast \geq 0 > s^{0\ast}$. Therefore, the proposition holds for time $t+n-1$. Next, we prove that the proposition holds for time $t+n-2, t+n-3, \cdots, t$. 

    Consider the state $s_{i^0(t+n-1)-1, t+n-2}$. Since $y^\ast(s_{i^0(t+n-1)-1, t+n-1}) = e_R$ and $y^\ast(s_{i^0(t+n-1), t+n-1}) = e_L$, by Lemma \ref{Lemmas0Difference}, if $s_{i^0(t+n-1)-1, t+n-2} \geq s^{0\ast}$ (the first case), then $y^\ast(s_{i^0(t+n-1)-1, t+n-2}) = e_L$; if $s_{i^0(t+n-1)-1, t+n-2} < s^{0\ast}$ (the second case),  then $y^\ast(s_{i^0(t+n-1)-1, t+n-2}) = e_R$. 
    
    For the first case, since $s_{i^0(t+n-1)-1, t+k} = s_{i^0(t+n-1)-1, t+n-2} + (n-2-k)\Delta_1$, we have $s_{i^0(t+n-1)-1, t+k} \geq s^{0\ast}$ for all $0\leq k\leq n-2$. For $0\leq k\leq n-2$, since $y^\ast(s_{i^0(t+n-1), t+k+1}) = y^\ast(s_{t+k+1}^0) = e_L$, we have $y^\ast(s_{i^0(t+n-1)-1, t+k}) = e_L$ by Lemma \ref{Lemmas0Difference}. Note that when $k=0$, $s_{i^0(t+n-1)-1, t} = s_t^0 = s_t^\ast$. Then the state sequence $s_t^\ast, s_{t+1}^\ast, \cdots, s_{t+n-1}^\ast$ equals to $s_{i^0(t+n-1)-1, t}, s_{i^0(t+n-1)-1, t+1}, \cdots, s_{i^0(t+n-1)-1, t+n-1}$. For $0\leq k\leq n-2$, we have $y^\ast(s_{t+k}^\ast) = y^\ast(s_{i^0(t+n-1)-1, t+k}) = e_L$ and $s_{t+k}^\ast = s_{i^0(t+n-1)-1, t+k} \geq s^{0\ast}$. Hence, the proposition holds in this case. 
    
    % For time $t+n-2$, we have $y^\ast(s_{i^0(t+n-1)-1, t+n-1}) = e_R$ and $s_{i^0(t+n-1)-1, t+n-1} = s_{t+n}^0 -\Delta_2 < -\Delta_2 < s^{0\ast}$. 

    For the second case, $y^\ast(s_{i^0(t+n-1)-1, t+n-2}) = e_R$. (i) If $n = 2$, then $y^\ast(s_t^\ast) = y^\ast(s_t^0) = e_R$ and $s_t^\ast = s_t^0 = s_{i^0(t+n-1)-1, t+n-2} < s^{0\ast}$. Hence, the proposition holds for interval $[t, t+n)$. (ii) If $n > 2$, then $s_{t+n-2}^0\geq 0$ and by Lemma \ref{LemmaProperties}, we have $i^0(t+n-1) = i^0(t+n-2)$ and $s_{i^0(t+n-1)-1, t+n-2} = s_{i^0(t+n-2)-1, t+n-2}$. Then, we already have $y^\ast(s_{i^0(t+n-2)-1, t+n-2}) = e_R$ and by Lemma \ref{Lemma3}, $y^\ast(s_{i^0(t+n-2), t+n-2}) = y^\ast(s_{t+n-2}^0) = e_L$. Thus, we are back to the same situation with time $t+n-1$, i.e., at this time, the optimal action at the state on the zero path is $e_L$ and the optimal action at the state which is at the left side of the zero path is $e_R$. Hence, the proposition holds for time $t+n-2$. Next, We prove that the proposition holds for time $t+n-3, t+n-4,  \cdots, t$.

    Consider the state $s_{i^0(t+n-1)-1, t+n-3}$. (i) If $s_{i^0(t+n-1)-1, t+n-3} \geq s^{0\ast}$, then we have $s_{i^0(t+n-1)-1, t+k} \geq s^{0\ast}$ for all $0\leq k\leq n-3$. The deduction is similar with that for the first case at time $t+n-2$ above. And the proposition holds. (ii) If $s_{i^0(t+n-1)-1, t+n-3} < s^{0\ast}$, the deduction is similar with that for the second case at time $t+n-2$ above and the proposition holds for time $t+n-3$. This process continues until time $t$.
    
    For time $t$, we consider the most special case in which $s_{i^0(t+n-1)-1, t+k} < s^{0\ast}$ for all $0\leq k\leq n-2$. From the above discussion of the second case, we have $y^\ast(s_{i^0(t+n-1)-1, t+k}) = e_R$ for all $0\leq k\leq n-2$. Note that when $k=0$, $s_{i^0(t+n-1)-1, t} = s_t^0 = s_t^\ast$. Then, the state $s_{t+1}^\ast = s_{t+1}^0$. Since $s_{t+k}^0\geq 0$ for $1\leq k < n$, by Lemma \ref{Lemma3} repeatedly, the state sequence $s_{t+1}^\ast, s_{t+2}^\ast, \cdots, s_{t+n-1}^\ast$ is $s_{t+1}^0, s_{t+2}^0, \cdots, s_{t+n-1}^0$. For $k=0$, we have $y^\ast(s_t^\ast) = y^\ast(s_{i^0(t+n-1)-1, t})  = e_R$ and $s_{i^0(t+n-1)-1, t} < s^{0\ast}$. For $1\leq k\leq n-2$, we have $y^\ast(s_{t+k}^\ast) = y^\ast(s_{i^0(t+n-1), t+k}) = y^\ast(s_{t+k}^0)  = e_L$ and $s_{t+k}^\ast \geq s^{0\ast}$. Hence, the proposition holds in this case.

    From the above discussion, we know that proposition holds from time $t$ to time $t+n$. By Lemma \ref{LemmaProperties} and \ref{Lemma4}, the states with the negative value on the zero path will appear periodically and they are also on the optimal path. The interval $[0, t^d)$ can be connected by such intervals like $[t, t+n)$. Therefore, for the state sequence $s_1^\ast, s_2^\ast, \cdots, s_{t^d}^\ast$ the proposition holds. 

\end{proof}

Below, we give two examples to illustrate the results of Theorem \ref{ThmBest}. 

\begin{figure}[htbp]
    \centering
    \includegraphics[width = 0.9\textwidth]{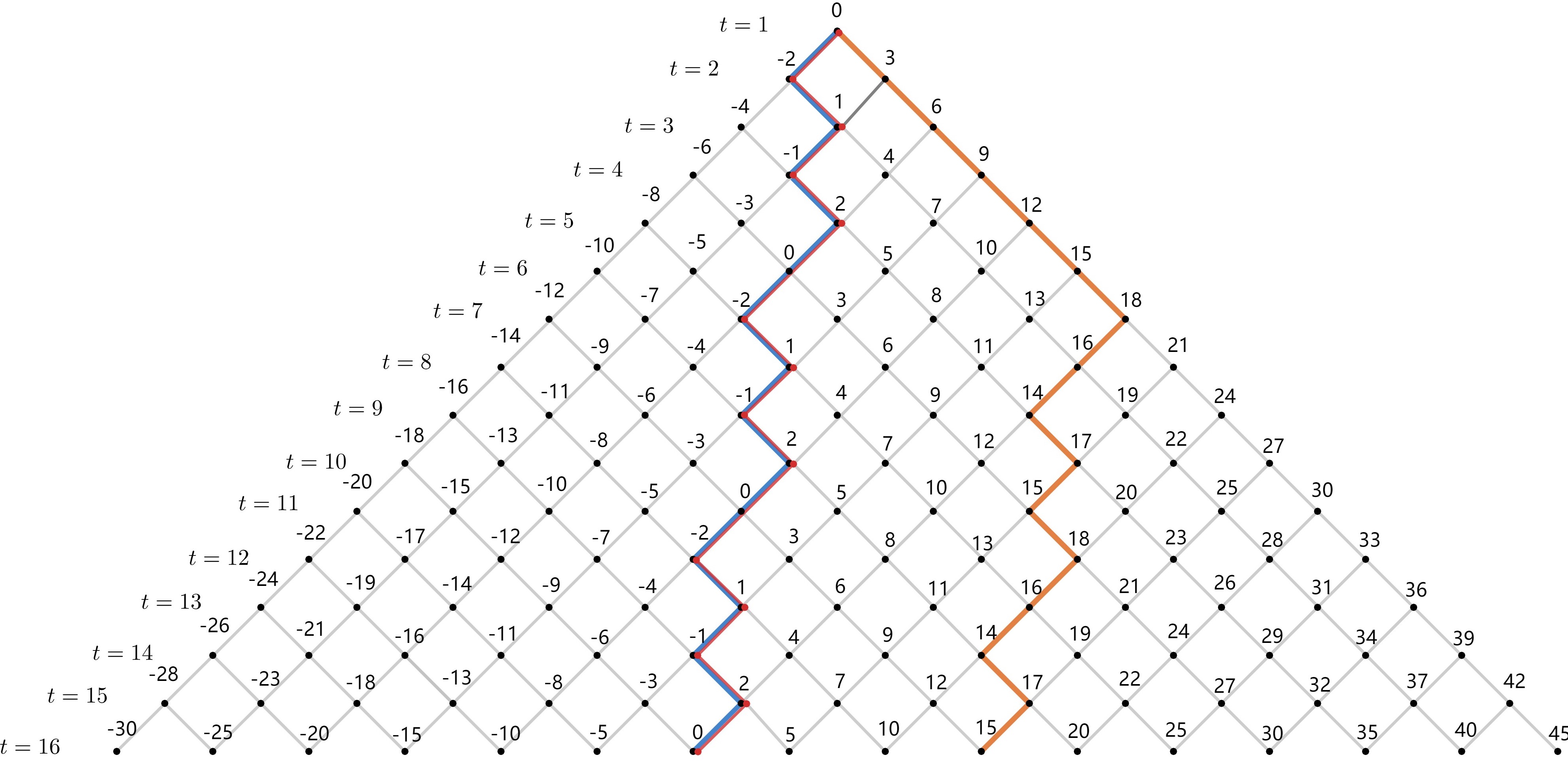}
    \caption{Illustration for Example \ref{Ex2}. In the graph, the orange polyline represents the myopic path, the blue polyline represents the zero path and the red polyline represents the optimal path. }
    \label{fig:Ex2}
\end{figure}

\begin{example}\label{Ex2} (Example \ref{EX1} continued.)
Take
$T=700$ and \begin{equation*}
    A = \begin{pmatrix}
    1 & 0\\
    -1 & 3
    \end{pmatrix}.
\end{equation*}
By simple calculation, 
$\Delta_1 = 2,\Delta_2=-3,\delta_1 = 1,\delta_2 = -4,s^\ast \approx 15.58$. 
The STTG as well as the myopic path, the zero path and the optimal path of this example are presented in Figure \ref{fig:Ex2}. 
From this figure, we can see that the optimal path is quite different from the myopic path. However, in this special setting that $\Delta_1 = 2,\Delta_2 = -3$, the optimal path is the same as the zero path from time $t=1$ to $t=t^d$.

Considering the payoffs, the time-averaged expected payoff of the optimal play is approximately 0.6667, which is bigger than the game value 0.6.
\end{example}

\begin{figure}[htbp]
    \centering
    \includegraphics[width = 0.9\textwidth]{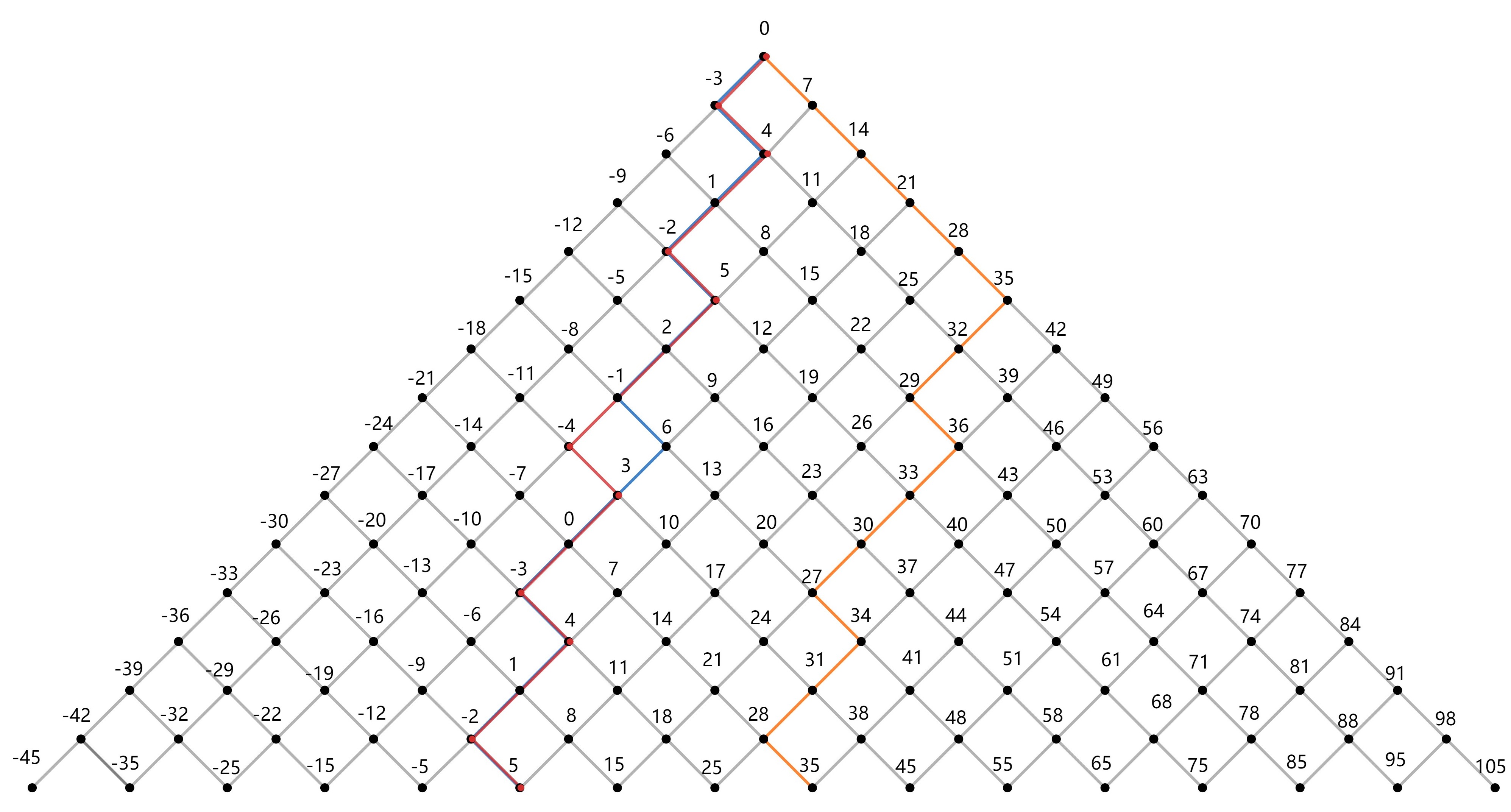}
    \caption{Illustration for Example \ref{EX3}. In the graph, the orange polyline represents the myopic path, the blue polyline represents the zero path and the red polyline represents the optimal path. }
    \label{fig:Ex3}
\end{figure}

\begin{example}\label{EX3}
Given that
$T=1000$ and \begin{equation*}
    A = \begin{pmatrix}
    1 & 0\\
    -2 & 7
    \end{pmatrix}.
\end{equation*}
By simple calculation, 
$\Delta_1 = 3,\Delta_2 = -7,\delta_1 = 1,\delta_2 = -9,s^\ast \approx 29.51, s^{0\ast} \approx -1.33.$ 
The STTG as well as the myopic path, the zero path and the optimal path of this example are presented in Figure \ref{fig:Ex3}. From Figure \ref{fig:Ex3}, we can see that the optimal path is quite different from the myopic path and the zero path. However, the optimal path has a lot of common states with the zero path. 

Considering the payoffs, the time-averaged expected payoff of the optimal play is approximately 0.8941, which is bigger than the game value 0.7.
\end{example}

Now, we have proven the property of the game system under the optimal play of player Y. 
Moreover, what is the time-averaged expected payoff of player Y under the optimal play?
What is the relation between it and the game value? 
Unfortunately, because of the exponential form of the stage strategy, we find it impossible to simplify the expression of the cumulative expected payoff. However, from the above examples, we observe that the time averaged payoff of player Y under the optimal play is bigger than the game value. In fact, this is generally valid and is stated in the proposition below.

\begin{proposition}
    The averaged payoff of player Y over the period is higher than the game value.
\end{proposition}

\begin{proof}
    For the payoff function $r(s, e_L)$ and $r(s, e_R)$, we find that
    \begin{equation}\label{formula_r_s_r}
        r(s, e_L) = \frac{a_{11}+a_{21}e^{-\eta s_t}}{e^{-\eta s_t}+1} = a_{21} + \frac{\Delta_1}{e^{-\eta s_t}+1}
    \end{equation}
    increases monotonically with $s$, while
    \begin{equation}\label{formula_r_s_l}
        r(s, e_R) = \frac{a_{12}+a_{22}e^{-\eta s_t}}{e^{-\eta s_t}+1} = a_{22} + \frac{\Delta_2}{e^{-\eta s_t}+1}
    \end{equation}
    decreases monotonically with $s$. This yields
    \begin{gather*}
        r(s, e_L)\geq r(s^{0\ast}, e_L), \quad \forall s\geq s^{0\ast},\\
        r(s, e_R) > r(s^{0\ast}, e_R),\quad \forall s < s^{0\ast}.
    \end{gather*}  
    By Proposition \ref{lemma_eL_state_greater_than_s0_er_state_less_than_el}, we know that the state with optimal action $e_R$ must be less than $s^{0\ast}$, and the state with optimal action $e_L$ must be not less than $s^{0\ast}$. In each period, the ratio of action $e_R$ to action $e_L$ is $\vert \Delta_1\vert : \vert \Delta_2\vert$. Hence, for the averaged payoff $\Bar{r}$ over the period, we obtain
    \begin{equation} \label{payoff_greater_than_value}
        \begin{aligned}
            \Bar{r} & = \frac{\sum_{s_t: y^\ast(s_t) = e_L} r(s_t, e_L) + \sum_{s_t: y^\ast(s_t) = e_R} r(s_t, e_R)}{T^\ast} \\ 
            & \geq \frac{\vert \Delta_1\vert r(s^{0\ast}, e_R)+\vert \Delta_2\vert r(s^{0\ast}, e_L)}{\vert \Delta_1\vert+\vert \Delta_2\vert} \\
            % & = \frac{a_{11}a_{22}-a_{12}a_{21}}{\vert \Delta_1\vert+\vert \Delta_2\vert}\\
            & = \frac{a_{11}a_{22}-a_{12}a_{21}}{a_{11}+a_{22}-a_{12}-a_{21}} \\
            & = v
        \end{aligned}
    \end{equation}
    where $T^\ast$ is the least period and $v$ is the game value. The second equality is obtained by substituting \eqref{formula_r_s_l} and \eqref{formula_r_s_r} into it.
\end{proof}

Note that at some stages in the period, the payoff of player Y is lower than the game value, but the averaged payoff over the period is higher. This implies that it is worth sacrificing the instantaneous payoffs at some stages for higher instantaneous payoffs at other stages. In fact, in many cases, the averaged payoff over the period is usually strictly greater than the game value.

\subsection{Extension to the Game With Rational Loss Values}\label{subsection_loss_matrix_elements_rational}
In this subsection, we will prove that Theorem \ref{ThmBest} also holds when the loss values are all rational numbers. Then, the loss matrix can be rewritten to be
\begin{equation*}
    A = 
    \begin{pmatrix}
    \frac{\tilde{a}_{11}}{q} & \frac{\tilde{a}_{12}}{q}\\
    \frac{\tilde{a}_{21}}{q} & \frac{\tilde{a}_{22}}{q}
    \end{pmatrix}
\end{equation*}
where $\tilde{a}_{11}, \tilde{a}_{12}, \tilde{a}_{21}, \tilde{a}_{22}, q$ are all integers.

Then, by \eqref{X_2divideX_1}, we have {
\begin{align*}
    \frac{x_{2,t}}{x_{1,t}} &= \frac{\exp (-\eta \sum_{\tau = 1}^{t-1}(\tilde{a}_{21}/q,\tilde{a}_{22}/q) y_\tau)}{\exp (-\eta \sum_{\tau = 1}^{t-1}(\tilde{a}_{11}/q,\tilde{a}_{12}/q) y_\tau)}\nonumber  \\
    &= \exp(-\frac{\eta}{q} (-\Delta_1,-\Delta_2) \sum_{\tau = 1}^{t-1} y_\tau)
\end{align*}
where $\tilde{\Delta}_1 = \tilde{a}_{11}-\tilde{a}_{21},\tilde{\Delta}_2 = \tilde{a}_{12}-\tilde{a}_{22}$. }

Let $\tilde{\eta} = \frac{\eta}{q}$ and $$\tilde{A} = \begin{pmatrix}
\tilde{a}_{11} & \tilde{a}_{12} \\ 
\tilde{a}_{21} & \tilde{a}_{22}.
\end{pmatrix}$$

Then, the optimal play problem of player Y against the Hedge algorithm with parameter $\eta$ in the game $A$ is equivalent to the optimal play problem of player Y against the Hedge algorithm with parameter {$\tilde{\eta}$} in the game {$\tilde{A}$}. Since $\eta = o(T)$,  we have $\tilde{\eta} = \frac{\eta}{q} = o(T)$ and thus the optimal action sequence of player Y is periodic over the whole horizon with a tiny interval truncated from the end, which means that Theorem \ref{ThmBest} still holds.

\subsection{Procedure to Obtain the Optimal Play}
Based on the main results, we give a simple method to obtain the optimal play of player Y. For simplicity, we still assume that the loss values are all integers. 

Given time $T$ and loss matrix $A$, we can calculate $\Delta_1 = a_{11}-a_{21}$, $\Delta_2=a_{12}-a_{22}$ and $T^\ast = \frac{m(\vert\Delta_1\vert,\vert\Delta_2\vert)}{\vert\Delta_1\vert}+\frac{m(\vert\Delta_1\vert,\vert\Delta_2\vert)}{\vert\Delta_2\vert}$. 

Then, we can compute the zero sequence $s_1^0,s_2^0,\cdots,s_{T+1}^0$. In fact, we only need to compute $s_1^0,s_2^0,\cdots, s_{T^\ast}^0$. 

Since $s_2^0<0$, we know $s_2^\ast = s_2^0$. Because of the periodicity of the optimal play, we solve 
\begin{equation*}
    \max\limits_{y_2,y_3,\cdots,y_{1+T^\ast}} f^\ast(s_2^\ast,s_{2+T^\ast}^\ast).
\end{equation*} 
to get the optimal play $y_2^\ast,y_3^\ast,\cdots,y_{1+T^\ast}^\ast$ and the sequence $s_2^\ast,s_3^\ast,\cdots, s_{2+T^\ast}^\ast$.

For the following periods, the optimal play is the same with the state sequence $y_2^\ast,y_3^\ast,\cdots,y_{1+T^\ast}^\ast$. 

Then, the periodic behavior of the optimal play ends at time $t^d$ and by \eqref{Rangetd}, time $t^d$ can be chosen as $t^d = T-T^\ast-\lceil\frac{s^\ast}{\vert\Delta_2\vert}\rceil$. According to the periodicity, we can compute $s_{t^d}^\ast$. 

Now, we only need to solve the Longest Path Problem starting from the state $s_{t^d}^\ast$ to $T+1$ stage. It is easy to solve because $T-t^d$ is quite small.

Compared with solving the Bellman Optimality Equation for all states, this method is simpler and more economical. Next, we give an example to show this method.

\begin{example}
Given $T=10000$ and \begin{equation*}
    A = \begin{pmatrix}
    1 & 0\\
    -1 & 3
    \end{pmatrix}.
\end{equation*}
By calculation, 
$\Delta_1 = 2,\Delta_2=-3,\delta_1 = 1,\delta_2 = -4,s^\ast \approx 58.87.$ The zero sequence is 0,-2,1,-1,2,0,-2,1,-1,2,0 $\cdots$.

Consider states at time $T$, state $s_{j^\ast,T}=57$ is the state such that $s_{j^\ast,T}<s^\ast$ and $s_{j^\ast+1,T}\geq s^\ast$. Then, consider $s_{j^\ast-T+t,t} = 57-3\times(T-t)$ and the zero sequence, we get $t^{\times} = T-19$ and $t^d = T-21$.
In this case, $T-t^d = 21$, which is quite small compared with $T=10000$. 
\end{example}

\section{{Some Discussions on Assumptions}}\label{sec6}

In the previous sections, we solved the problem proposed in Section \ref{sec2} based on some assumptions. 
In this section, we will discuss the impact of these assumptions and see what would happen if some assumptions are not satisfied.

\subsection{{If the Loss Matrix Elements Are Irrational Numbers}}\label{subsection_irrational_numbers}

In the above sections, we prove the periodicity of the optimal action sequence of player Y against the Hedge algorithm, which helps us find it efficiently. 
However, this nice property may not hold for the game with irrational loss values. In fact, the parameter $\Delta_2/\Delta_1$ is essential to decide whether periodicity holds.

Why? The reason lies in the updating formula of the state $s_t$. By \eqref{StateTransFormula}, if player Y takes action $e_L$, then $s_{t+1} = s_t - \Delta_1$, and if player Y takes action $e_R$, then $s_{t+1} = s_t - \Delta_2$.
Since $\Delta_1\Delta_2<0$, if $\Delta_1/\Delta_2$ is rational, then there exist two integers $p$ and $q$ such that $p \Delta_1 + q\Delta_2 = 0$, which implies that after $p+q$ stages, the state will be the same with the original state. However, if $\Delta_2/\Delta_1$ is irrational, then for any two integers $p$ and $q$, $p \Delta_1 + q\Delta_2$ will never be zero, implying that the state will never come back to the original state. Thus, periodicity does not hold for any action sequence of player Y.

So, how can we find the optimal play of player Y if $\Delta_1/\Delta_2$ is irrational? A natural idea is to use the optimal play for rational numbers to approximate the optimal play for irrational numbers. This needs to get a more explicit description of the optimal strategy for rational numbers. The zero point $s^{0\ast}$ of function \eqref{f(s)Def} is crucial, around which the optimal action sequence oscillates. We guess that this still holds for irrational $\Delta_1/\Delta_2$ and the optimal play can be constructed from this phenomenon.

\subsection{{If Player X Has a Dominant Action}}\label{subsection_player_x_has_dominant_strategy}

In the above sections, we assumed that player X has no dominant action, i.e., $\Delta_1\Delta_2<0$. In this subsection, we will prove that if player X has a dominant action, then no matter whether player Y has the dominant action, we can still find the optimal action sequence for player Y. 
We state this in Proposition \ref{Optimal>0}.

\begin{proposition}\label{Optimal>0}

If $\Delta_1>0,\Delta_2>0$, i.e., action D is the dominant action of player X, then there exists time $t^m$ such that $y_{t}^\ast=e_R$ for all time $t\geq t^m$.
\end{proposition}

\begin{proof}

\textbf{Case 1}: if player Y has a dominant action, i.e., $\delta_1\delta_2>0$. If player Y has a dominant action, taking the dominated action will not improve the expected payoff of taking the dominant action. When player X has a dominant action, player Y taking action L and R will both make player X tend to action $D$. Then, the expected payoff of taking the dominant action for player Y decreases. Therefore, the optimal choice for player Y is to always take the dominant action. Take $t^m = 1$, and the proposition holds in this case.

\textbf{Case 2}: if player Y does not have a dominant action, i.e., $\delta_1\delta_2<0$. 
Recall that $s^\ast = -\ln \left(-{\delta_1}/{\delta_2}\right)/\eta.$ 
By \eqref{PayoffFunc}, we have
\begin{equation*}
    r(s,e_R) = \frac{a_{12}+e^{-\eta s}a_{22}}{1+e^{-\eta s}} = a_{22}+\frac{\Delta_2}{1+e^{-\eta s}}.
\end{equation*}
and
$$r(s,e_R)-r(s,e_L) = -\frac{\delta_1+e^{-\eta s}\delta_2}{1+e^{-\eta s}}.$$ 
Then, we have $r(s,e_R)$ increases with $s$, and
\begin{equation}\label{Toprove1}
    \left\{ \begin{array}{cc}
    r(s,e_R)>r(s,e_L),     &\quad \text{if}\ s<s^\ast;  \\
    r(s,e_R)\leq r(s,e_L),     & \quad \text{if}\ s\geq s^\ast.    
    \end{array}
\right.
\end{equation}

Firstly, we prove that there exists time $t^m$ such that for all $t\geq t^m$ and $1\leq i\leq t,$
$s_{i,t}\leq s^\ast. $ We prove this in two cases.
(1) If $s^\ast> 0$, $s_{i,t} = -(t-i)\Delta_1-(i-1)\Delta_2<0$ for all $i$ and $t$. Obviously, such $t^m$ exists. (2)
If $s^\ast\leq 0$, since $s_{t,t} = -(t-1)\Delta_2$, there exists time $t^m$ such that $s_{t^m,t^m}\leq s^\ast$. Further, for all $t\geq t^m$, we have $s_{i,t}\leq s_{t,t}\leq s_{t^m,t^m}$ for all $1\leq i\leq t$. 
Therefore, in both cases, there exists time $t^m$ such that for all $t\geq t^m$ 
$s_{i,t}\leq s^\ast $ for all $1\leq i\leq t$. 

Secondly, we prove that for all $t\geq t^m$ and $1\leq i\leq t$, $y^\ast(s_{i,t}) = e_R$ by mathematical induction. 

First, we prove that for $t=T$, the above claim holds.
Since $T\geq t^m$, we have just proved that $s_{i,T}\leq s^\ast$ for all $1\leq i\leq T$. Since 
$f^\ast(s_{T+1})= 0$ for all $s_{T+1}$, by \eqref{Toprove1}, we have $y^\ast(s_{i,T})=e_R$ for all $1\leq i\leq T$.

Second, assume that for $t$ such that $t^m+1\leq t\leq T$, we have $y^\ast(s_{i,t})=e_R$ for all $1\leq i\leq t$. Then, we need to prove that $y^\ast(s_{i,t-1})=e_R$.
Since $r(s,e_R)$ increases with $s$ and $s_{i,t-1}\leq s^\ast$, by \eqref{Toprove1}, we have
\begin{equation*}
    \begin{aligned}
    r(s_{i,t-1},e_R)+f^\ast(h(s_{i,t-1},e_R)) &=  r(s_{i,t-1},e_R)+f^\ast(s_{i+1,t})\\
    & = r(s_{i,t-1},e_R) + \sum_{\tau = t}^T r(s_{i-t+\tau+1,\tau},e_R)\\
    &> r(s_{i,t-1},e_L) + \sum_{\tau = t}^T r(s_{i-t+\tau,\tau},e_R)\\
    & = r(s_{i,t-1},e_L)+f^\ast(s_{i,t})\\
    & = r(s_{i,t-1},e_L)+f^\ast(h(s_{i,t-1},e_L))
    \end{aligned}
\end{equation*}
where the second and the third equality hold because of the induction assumption.
Thus, we have $y^\ast(s_{i,t-1}) = e_R$ for all $1\leq i\leq t-1$.

Therefore, for all $t\geq t^m$, we have $y^\ast(s_{i,t})=e_R$ for all $1\leq i\leq t$. Then,  $y_{t}^\ast=e_R$ for all time $t\geq t^m$.

\end{proof}

\subsection{{If Player Y Has a Dominant Action}}\label{subsection_player_Y_has_dominant_strategy}

In the above subsection, we proved the optimal action sequence of player Y when player X has a dominant action. In this subsection, we will show that when player Y has a dominant action, the situation is very different and makes the problem much more complicated.

What makes the difference? In order to get the main results, we tried to solve the Bellman Optimality Equation by backward induction. If player Y does not have a dominant action, there exists some index $i$ such that $ y^\ast(s_{j, T}) = e_L$ for all indexes $j \geq i$ and $ y^\ast(s_{j, T}) = e_R$ for all indexes $j < i$. Then, we can use Corollary \ref{LemmaRtoRLtoL} to recursively solve the Bellman Optimality Equation from bottom to top on STTG.  However, if player Y has a dominant action, all the solutions at time $T$ are the dominant action. Then, we can only use part of Corollary \ref{LemmaRtoRLtoL} and the optimal action at more than half of states on STTG cannot be obtained by recursion, which makes it difficult to locate the optimal path.

So, what is the optimal play for player Y now? In fact, the problem is very complicated, which can be shown by the examples below.

\begin{figure}[htbp]
    \centering
    \includegraphics[width = 0.9\textwidth]{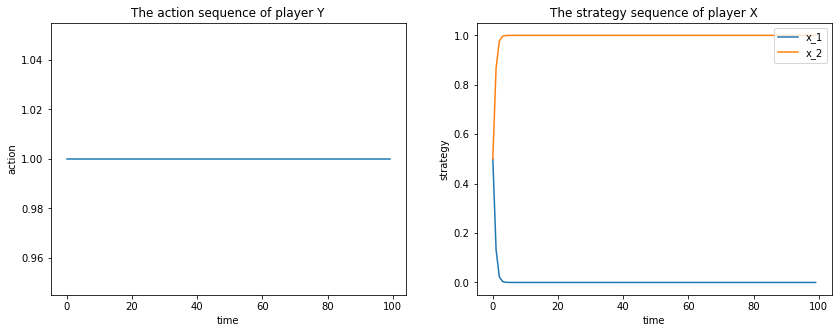}
    \caption{ {An example when always taking the dominant action is optimal for player Y in repeated game.}}
    \label{fig:dominant_strategy_example_1}
\end{figure}

\begin{example}\label{fig_dominant_strategy_is_optimal}

    Take $T = 100$ and the loss matrix is
    \begin{equation*}
    A = \begin{pmatrix}
    0 & 10\\
    1 & 2
    \end{pmatrix}.
\end{equation*}
Here, player Y has a dominant action $R$. In this example, by the Bellman Optimality Equation, the optimal choice for player Y is to always take the dominant action. The optimal action sequence of player Y and the strategy sequence of player X are shown in Figure \ref{fig:dominant_strategy_example_1}, in which $y = 1$ represents the action $R$ of player Y. 
\end{example}

\begin{figure}[htbp]
    \centering
    \includegraphics[width = 0.9\textwidth]{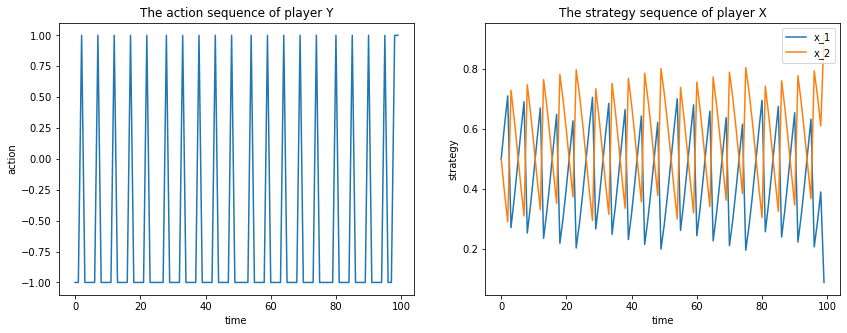}
    \caption{ {An example when always taking the dominant action is not optimal for player Y in repeated game.}}
    \label{fig:dominant_strategy_example_2}
\end{figure}

\begin{example}\label{fig_dominant_is_not_optimal_rational}

Take $T = 100$ and the loss matrix is
    \begin{equation*}
    A = \begin{pmatrix}
    0 & 10\\
    1.9 & 2
    \end{pmatrix},
\end{equation*}
whose elements are the same as Example \ref{fig_dominant_strategy_is_optimal} except $a_{21}$. However, in this example, always taking the dominant action is not the optimal choice for player Y!  The optimal action sequence of player Y and the strategy sequence of player X are shown in Figure \ref{fig:dominant_strategy_example_2}, in which $y = 1$ represents the action $R$ of player Y and $y = -1$ represents the action $L$. 

\end{example}

\begin{figure}[htbp]
    \centering
    \includegraphics[width = 0.9\textwidth]{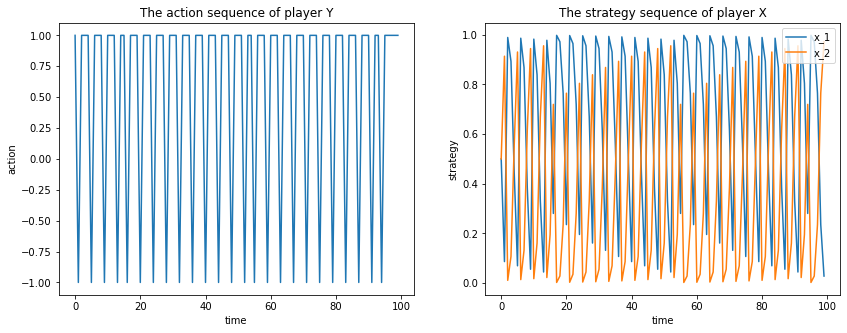}
    \caption{ {An example when always taking the dominant action is not optimal for player Y and the loss matrix elements are all integers.}}
    \label{fig:dominant_strategy_example_3}
\end{figure}

\begin{example}
Take $T = 100$ and the loss matrix is
    \begin{equation*}
    A = \begin{pmatrix}
    0 & 30\\
    19 & 20
    \end{pmatrix}.
\end{equation*}
In this example, always taking the dominant action either is not the optimal choice for player Y. The optimal action sequence and the strategy sequence of player X are shown in Figure \ref{fig:dominant_strategy_example_3}. 
This example shows that even when the elements of the loss matrix are all integers, it is still possible that always taking the dominant action is not optimal for player Y.
\end{example}

From the above examples, we can see that the optimal action sequence of player Y depends on the elements of the loss matrix in a sensitive way. For instance, the loss matrix in Example \ref{fig_dominant_is_not_optimal_rational} differs from that in Example \ref{fig_dominant_strategy_is_optimal} in only one element, but the optimal action sequences are very different.

So, what is the optimal play for player Y? The essential problem for player Y is to tackle a trade-off between the instantaneous and long-run payoffs. Because of the exponential form of the stage strategy, it is very hard to compare the cumulative expected payoff on different paths. Therefore, in order to get the optimal path, more detailed analysis and comparison are necessary and we leave it as future work.

\subsection{{The Effect of the Initial Strategy}}\label{subsection_effect_of_initial_strategy}

In this paper, we assume that the initial strategy $x_1$ of player X is $(\frac12, \frac12)$. This is a common setting for the Hedge algorithm. Below, we will illustrate that if the initial strategy changes, the optimal action sequence of player Y is still the same except for the first few actions.

If the initial strategy is not $(\frac12, \frac12)$, the IEP $x_1^T A y_1$ of the first stage is different, and thus the lengths of the two edges connecting the first node $s_{1}$ are different too. However, from the updating formula \eqref{HedgeStrategyFormula} of the Hedge algorithm, the stage strategy $x_t$ of player X is entirely determined by the action sequence $ y_1, y_2, \cdots, y_{t-1}$ of player Y, so $x_1$ does not affect the strategy of player X at time $t\geq 2$. Hence, if the initial strategy is not $(\frac12, \frac12)$, the corresponding states of nodes and lengths of the edges from the 2nd to (T+1)th row on STTG are still the same as well as the Bellman Optimality Equation and its solution (i.e., the optimal action \eqref{yast}) for each state.

By Lemma \ref{Lemma4}, for $t>2$, if $s_t^0 < 0$ , then $s_t^\ast = s_t^0$. That is to say, if the state at time $t > 2$ on the zero path is negative, then we know this state will be the state at time $t$ on the new optimal path. Since the zero path remains the same as that under the new initial strategy, no matter which action player Y takes at the first stage, the state on the zero path at time $t$ satisfying $s_t^0<0$ is still on the new optimal path. Therefore, $x_1$ can affect only the first few actions of the optimal action sequence.

\section{Conclusion and Future Work}\label{sec7}
In this paper, we study the optimal play for one player in a two-player finitely-repeated $2\times 2$ zero-sum game, given that the other player employs the Hedge algorithm to update its stage strategy.            
We propose the State Transition Triangle Graph and find that the problem of the optimal play against Hedge is equivalent to solving the Longest Path Problem on STTG. Then, we study the property of the game system under the myopic and optimal play. First, we prove that the myopic play will enter a cycle after small stages. Specifically, the zero path shows periodic behavior throughout the whole time horizon. Then, we get the recurrence relation between the optimal actions of states at adjacent times. Based on these findings, if no dominant action exists, we prove that the optimal play is periodic on the time interval truncated by a tiny segment at the end.

The optimal play against the adaptive algorithm in repeated games will become much more significant and there are still many open problems, which are worthy of being investigated. For example, if the Hedge algorithm has a decaying learning rate, which seems a reasonable assumption to avoid the prior knowledge about the game, we will encounter the time-varying system and the analysis will become very difficult. 
From our recent work, we guess that similar results would hold for general games with more actions under some conditions on the loss values, only that the optimal path would be more complicated. Such extensions are attractive and also challenging. Besides, if the specific forms of the adaptive algorithms are not known in advance, it seems possible to identify them and such identification problems for learning algorithms seem very rare in the literature. All these problems need our deeper understanding of the game system driven by the learning algorithms. 
We leave them as future work.

\section{Acknowledgement}

This work was supported by the Strategic Priority Research Program of Chinese Academy
of Sciences under Grant No. XDA27000000, the Natural Science Foundation of China under Grant T2293770, the National Key Research and Development Program of China under grant No.2022YFA1004600, the Major Project on New Generation of Artificial Intelligence from the Ministry of Science and Technology (MOST) of China under Grant No. 2018AAA0101002.

%Bibliography
\bibliographystyle{unsrt}  
\bibliography{ref}

\end{document}